\documentclass[english,12pt]{bourbaki}
\usepackage{picinpar}

\usepackage[latin1]{inputenc}

1

\newtheorem{lem}{Lemma}[section]
\newtheorem{theor}[lem]{Theorem}
\newtheorem{corol}[lem]{Corollary}
\newtheorem{propo}[lem]{Proposition}
\newtheorem{remar}[lem]{Remark}
\newtheorem{defin}[lem]{Definition}
\newtheorem{conje}[lem]{Conjecture}

\newenvironment{lemma}
{\begin{lem}\sl } {\end{lem}}

\newenvironment{theorem}
{\begin{theor}\sl} {\end{theor}}

\newenvironment{rtheorem}[1]
{\begin{theor}[{\rm #1}]\sl} {\end{theor}}

\newenvironment{corollary}
{\begin{corol}\sl } {\end{corol}}

\newenvironment{rcorollary}[1]
{\begin{corol}[{\rm #1}]\sl } {\end{corol}}

\newenvironment{proposition}
{\begin{propo}\sl } {\end{propo}}

\newenvironment{remark}
{\begin{remar} \rm \footnotesize} {\end{remar}}

\newcommand\eps\varepsilon
\newcommand\ph\varphi
\newcommand\A{{\mathbb A}}
\renewcommand\P{{\mathbb P}}
\newcommand\Q{{\mathbb Q}}
\newcommand\barQ{{\bar\Q}}
\newcommand\R{{\mathbb R}}
\newcommand\C{{\mathbb C}}
\newcommand\G{{\mathbb G}}
\newcommand\Z{{\mathbb Z}}
\newcommand\ba{{\mathbf a}}
\newcommand\bx{{\mathbf x}}
\newcommand\by{{\mathbf y}}

\renewcommand\AA{{\mathcal A}}
\newcommand\OO{{\mathcal O}}
\newcommand\QQ{{\mathcal Q}}

\newcommand\LL{{\mathcal L}}

\newcommand\ord{{\mathrm {ord}}}

\newcommand\genus{{\mathbf g}}

\newcommand\gerp{{\mathfrak p}}

\renewcommand\mod{{\rm \,mod\,}}

\newlength{\tmtwo}
\date{Novembre 2006}
\bbkannee{59\`eme ann\'ee, 2006-2007}
\bbknumero{967}\title{The Many Faces of the Subspace Theorem}
\subtitle{after Adamczewski, Bugeaud, Corvaja, Zannier\dots}
\author{{Yuri F. BILU}}
\address{{Universit\'e Bordeaux I\\
UFR de Math\'ematiques et Informatique\\
A2X~:
Laboratoire de Théorie des Nombres\\
et d'Algorithmique Arithmétique\\
(UMR 5465 du CNRS)\\
351 cours de la Lib\'eration\\
F-33405 Talence Cedex}}\email{{yuri@math.u-bordeaux1.fr}}

\begin{document}

\maketitle

{\scriptsize \tableofcontents

\begin{flushright}
\parbox{0.5\textwidth}{\textit{And we discovered subspace. It gave us our galaxy and it gave us the universe. And we saw other advanced life. And we subdued it or we crushed it\dots With subspace, our empire would surely know no boundaries.}}

\smallskip\noindent
(From \textit{The Great War} computer game)
\end{flushright}

}

\section{Introduction}
This is not a typical Bourbaki talk. A generic \textit{exposé} on this seminar is, normally, a report on a recent seminal achievement, usually involving new technique. The principal character of this talk is the Subspace Theorem of Wolfgang Schmidt, known for almost forty years. All results I am going to talk about rely on this celebrated theorem (more precisely, on the generalization due to Hans Peter Schlickewei). Moreover, in all cases it is by far the most significant ingredient of the proof.

Of course, the last remark is not meant to belittle the work of the authors of the results I am going to speak about. Adapting the Subspace Theorem to a concrete problem is often a formidable task, requiring great imagination and ingenuity. 

During the last decade the Subspace Theorem found  several quite unexpected applications, mainly in the Diophantine Analysis and in the  Transcendence Theory. Among the great variety of spectacular results, I have chosen several which are technically simpler  and which allow one to appreciate how miraculously does the Subspace Theorem  emerge in numerous  situations, implying beautiful  solutions to difficult problems hardly anybody hoped to solve so easily. 

The three main topics discussed in this article are: 

\begin{itemize}

\item
the work of Adamczewski and Bugeaud on complexity of algebraic numbers;

\item
the work of Corvaja and Zannier on Diophantine equations with power sums;

\item 
the work of Corvaja and Zannier on integral  points on curves and surfaces, and the subsequent development due to Levin and Autissier. 

\end{itemize}

In particular, we give a complete proof of the  beautiful theorem of Levin and Autissier (see Theorem~\ref{tlev}): \textsl{an affine surface with~$4$ (or more) properly intersecting ample divisors at infinity cannot have a Zariski dense set of integral points. }

\begin{sloppypar}
	
	Originally, Schmidt proved his theorem for the needs of two important subjects: norm form equations and exponential Diophantine equations (including the polynomial-exponential equations and linear recurrence sequences). These ``traditional'' applications of the Subspace Theorem  form a vast subject, interesting on its own; we do not discuss it here (except for a few motivating remarks in Section~\ref{spowsum}).  Neither do we discuss the quantitative aspect of the Subspace Theorem. For this, the reader should consult the  fundamental work of Evertse and Schlickewei (see \cite{ES99,ES02,Sc96,Sc98,Sc03} and the references therein).  
\end{sloppypar}

Some of the results stated here admit far-going generalizations, but  I do not always mention them: the purpose of this talk is to exhibit ideas rather than  to survey the best known results. 

In Section~\ref{sst} we introduce the Subspace Theorem. Sections~\ref{scomp},~\ref{spowsum} and~\ref{sipon} are totally independent and can be read in any order.

\section{The Subspace Theorem}
\label{sst}
In this section we give a statement of the Subspace Theorem. Before formulating it in full generality, we consider several particular cases, to make the general case more motivated.

\subsection{The Theorem of Roth}
\label{ssroth}

In 1955 K.~F.~Roth~\cite{Ro55} proved that algebraic numbers cannot be ``well approximated'' by rationals.
\begin{rtheorem}{Roth}
Let~$\alpha$ be an irrational 	algebraic number. Then for any   ${\eps>0}$ the inequality 
$$
\left|\alpha-\frac yx\right|<\frac1{|x|^{2+\eps}}
$$
has only finitely many solutions in non-zero ${x,y\in \Z}$. 
\end{rtheorem}

This result is, in a sense, best possible, because, by the Dirichlet approximation theorem, the inequality ${|\alpha-y/x|\le |x|^{-2}}$ has infinitely many solutions.

The theorem of Roth has a glorious history. Already Liouville showed in 1844 the inequality ${|\alpha-y/x|\ge c(\alpha)|x|^{-n}}$, where~$n$ is the degree of the algebraic number~$\alpha$, and used this to give first examples of transcendental numbers. However, Liouville's theorem was too weak for serious applications in the Diophantine Analysis. In 1909 A.~Thue~\cite{Th09} made a breakthrough, proving that ${|\alpha-y/x|\le|x|^{-n/2-1-\eps}}$ has finitely many solutions. A series of refinements (most notable being due to Siegel~\cite{Si29}) followed, and Roth made the final (though very important and difficult) step.

Kurt Mahler, who was a long proponent of $p$-adic Diophantine approximations, suggested to his student D.~Ridout~\cite{Ri58} to extend Roth's theorem to the non-archimedean domain. To state Ridout's result, we need to introduce some notation. For every prime number~$p$, including the ``infinite prime'' ${p=\infty}$, we let ${|\cdot|_p}$ be the usual $p$-adic norm on~$\Q$ (so that ${|p|_p=p^{-1}}$ if ${p<\infty}$ and ${|2006|_\infty =2006}$), somehow extended to the algebraic closure~$\bar \Q$.  For a rational number ${\xi=y/x}$ with ${\gcd(x,y)=1}$ we define its \textit{height}  by 
\begin{equation}
\label{eher}
H(\xi)=\max\{|x|, |y|\}.
\end{equation}
One immediately verifies that 
\begin{equation}
\label{ehaffq}
H(\xi)=\prod_p\max\left\{1, |\xi|_p\right\}=\left(\prod_p\min\left\{1, |\xi|_p\right\}\right)^{-1},
\end{equation}
where the products extend to all prime numbers, including the infinite prime. 

Now let~$S$ be a finite set of primes, including ${p=\infty}$, and for every ${p\in S}$ we fix an  algebraic number~$\alpha_p$. Ridout proved that for any ${\eps>0}$ the inequality 
$$
\prod_{p\in S}\min\left\{1,\left|\alpha_p-\xi\right|_p\right\}<\frac1{H(\xi)^{2+\eps}}
$$
has finitely many solutions in ${\xi\in \Q}$.

While the theorem of Roth becomes interesting only when the degree of~$\alpha$ is at least~$3$, the theorem of Ridout is quite non-trivial even when the  ``targets''~$\alpha_p$ are rational. Moreover, one can also allow ``infinite'' targets, with the standard convention ${\infty -\xi=\xi^{-1}}$. The following particular case of Ridout's theorem is especially useful: given an algebraic number~$\alpha$, a set~$S$ of prime numbers, and ${\eps>0}$, the inequality 
$$
|\alpha-\xi|<H(\xi)^{-1-\eps}
$$
has finitely many solutions in $S$-integers\footnote{A rational number is called $S$-integer if its denominator is divisible only by the prime numbers from~$S$.}~$\xi$. To prove this, consider the theorem of Ridout with ${\alpha_\infty=\alpha}$ and with ${\alpha_p=\infty}$ for ${p\ne \infty}$, and apply~(\ref{ehaffq}).

One consequence of this result is that the decimal expansion of an algebraic number cannot have ``too long'' blocks of zeros. More precisely, let ${0.a_1a_2\ldots}$ be the decimal expansion of an algebraic number, and for every~$n$ define $\ell(n)$ as the minimal ${\ell\ge0}$ such that ${a_{n+\ell}\ne 0}$; then  ${\ell(n)=o(n)}$ as ${n\to \infty}$. To show this, apply the above-stated particular case of the theorem of Ridout with ${S=\{2,5,\infty\}}$. More generally, the decimal expansion of an algebraic number cannot have ``too long'' periodic blocks.

S.~Lang extended the theorem of Roth-Ridout to approximation of algebraic numbers by the elements of a given number field. We invite the reader to consult Chapter~7 of his book~\cite{La83} or Part~D of the more recent volume~\cite{HS00} for the statement and the proof of Lang's theorem.

\subsection{The Statement of the Subspace Theorem}

Now we have enough motivation to state the Subspace Theorem.
We begin with the original theorem of Schmidt~\cite{Sc72} (see also~\cite{Sc80} for a very detailed proof).

\begin{sloppypar}
	\begin{rtheorem}{W.~M.~Schmidt}
	Let ${L_1, \ldots , L_m}$ be linearly independent linear forms in~$m$ variables with (real) algebraic coefficients. Then for any ${\eps>0}$  the  solutions ${\bx=(x_1, \ldots, x_m) \in \Z^m}$ of the inequality
	$$
	\left|L_1(\bx) \cdots  L_m(\bx)\right|\le \|\bx\|^{-\eps}
	$$
	are contained in  finitely many proper linear subspaces of~$\Q^m$. (Here ${\|\bx\|=\max_i\{|x_i|\}}$.)
	\end{rtheorem}
\end{sloppypar}

Putting ${m=2}$, ${L_1(x,y)=x\alpha-y}$ and ${L_2(x,y)=x}$, we recover the theorem of Roth. 

The theorem of Schmidt is not sufficient for many applications. One needs a non-archimedean generalization of it, analogous to Ridout's generalization of Roth's theorem. This result was obtained by Schlickewei~\cite{Sc76,Sc76a}. 
As in the previous section, let~$S$ be a finite set of prime numbers, including ${p=\infty}$, and pick an extension of every $p$-adic valuation to~$\bar\Q$.

\begin{rtheorem}{H.~P.~Schlickewei}
\label{tschl}
For every ${p\in S}$ let  ${L_{1,p}, \ldots , L_{m,p}}$ be linearly independent linear forms in~$m$ variables with algebraic coefficients.
Then for any ${\eps>0}$  the  solutions ${\bx \in \Z^m}$ of the inequality
$$
\prod_{p\in S}\prod_{i=1}^m\left|L_{i,p}(\bx)\right|_p\le \|\bx\|^{-\eps}
$$
are contained in  finitely many proper linear subspaces of~$\Q^m$. 
\end{rtheorem}

It is usually more convenient to allow the variables ${x_1, \ldots, x_m}$ to be $S$-integers rather than integers.  To restate Schlickewei's theorem using the $S$-integer variables, one needs an adequate measure of the ``size'' of a vector with $S$-integer (or, more generally, rational) coordinates; evidently,   the sup-norm ${\|\bx\|}$ cannot serve for this purpose. Thus, let~$\bx$ be a non-zero vector from~$\Q^m$; we define its \textit{height} by
\begin{equation}
\label{eheq}
H(\bx)= \prod_p \|\bx\|_p,
\end{equation}
where ${\|\bx\|_p = \max\{|x_1|_p, \ldots, |x_m|_p\}}$, 
and the product extends to all rational primes, including ${p=\infty}$. 

\begin{sloppypar}
	The height function, defined this way, is ``projective'': if ${a\in \Q^\ast}$ then ${H(a\bx)=H(\bx)}$ (this is an immediate consequence of the product formula). When the coordinates ${x_1, \ldots, x_m}$ are coprime integers, we have ${H(\bx)=\|\bx\|}$. 
\end{sloppypar}

\begin{remark}
One piece of warning: the height of a rational number~$\xi$, defined in~(\ref{eher}) is \emph{not} equal to the height of the ``one-dimensional vector'' with the coordinate~$\xi$; in fact, the height of a non-zero one-dimensional vector is~$1$, by the product formula,  while  $H(\xi)$ is  the height of the $2$-dimensional vector $(1, \xi)$, according to~(\ref{ehaffq}). This abuse of notation is quite common and will not lead to any confusion.
\end{remark}

Denote by~$\Z_S$ the ring of $S$-integers.
Now Theorem~\ref{tschl} can be re-stated as follows.

\paragraph*{\textsc{Theorem~\ref{tschl}${}'$}} \textsl{In the set-up of Theorem~\ref{tschl}, the  solutions ${\bx \in \Z_S^m}$ of the inequality
$$
\prod_{p\in S}\prod_{i=1}^m\left|L_{i,p}(\bx)\right|_p\le H(\bx)^{-\eps}
$$
are contained in  finitely many proper linear subspaces of~$\Q^m$.} 

It is very easy to deduce Theorem~\ref{tschl}${}'$ from Theorem~\ref{tschl}; we leave this as an exercise for the reader. (One should use the ``product formula'' ${\prod_p|a|_p=1}$, where ${a\in \Q^\ast}$ and the product extends to all rational primes, including ${p=\infty}$.)

Unfortunately, for many applications  Theorem~\ref{tschl}${}'$ is insufficient as well: one needs to extend it to the case when the variables ${x_1, \ldots, x_m}$ belong to an arbitrary number field. This was also done by Schlickewei~\cite{Sc77}. Before stating the theorem, we need to make some conventions. Let~$K$ be a number field of degree ${d=[K:\Q]}$  and let~$M_K$ be the set of all absolute values on~$K$. Recall that the set~$M_K$ consists of infinitely many \textit{finite} absolute values, corresponding to prime ideals of the field~$K$, and finitely many \textit{infinite} absolute values, corresponding to real embeddings of~$K$ (real absolute values) and pairs of complex conjugate embeddings (complex absolute values). 

We normalize the absolute values on~$K$ as follows.  If ${v \in M_K}$ is a $\gerp$-adic absolute value,  then we normalize it so that  ${|p|_v=p^{-d_v/d}}$, where~$p$ is the prime number below the prime ideal~$\gerp$ and  ${d_v=[K_v:\Q_p]}$ is the local degree. If~$v$ is an infinite absolute value, then we normalize it to have ${|2006|_v= 2006^{d_v/d}}$, where~$d_v$ is again the local degree (that is, ${d_v=1}$ if~$v$ is real and ${d_v=2}$ if~$v$ is complex). With this normalization we have the product formula in the form ${\prod_{v\in M_K}|a|_v=1}$, where ${a \in K^\ast}$.  

We also need to define the height of a vector ${\bx\in K^m}$. By analogy with~(\ref{eheq}) we put
${H(\bx)=\prod_{v\in M_K}\|\bx\|_v}$,
where ${\|\bx\|_v= \max\{|x_1|_v, \ldots, |x_m|_v\}}$. An easy verification shows that for ${\bx\in \Q^m}$ this definition agrees with~(\ref{eheq}).

Now we are ready to state the Subspace Theorem in its most general form. Let~$K$ be a number field, and let~$S$ be a finite set of absolute values of~$K$ (normalized as above), including all the infinite absolute values. We denote by~$\OO_S$ the ring of $S$-integers\footnote{An element ${\alpha \in K}$ is called $S$-integer if ${|\alpha|_v\le 1}$ for all ${v\notin S}$.} of the field~$K$. 

\begin{rtheorem}{H.~P.~Schlickewei}
\label{tssnf}
For every ${v\in S}$ let  ${L_{1,v}, \ldots , L_{m,v}}$ be linearly independent linear forms in~$m$ variables with algebraic coefficients.
Then for any ${\eps>0}$  the  solutions ${\bx \in \OO_S^m}$ of the inequality
$$
\prod_{v\in S}\prod_{i=1}^m\left|L_{i,v}(\bx)\right|_v\le H(\bx)^{-\eps}
$$
are contained in  finitely many proper linear subspaces of~$K^m$. 
\end{rtheorem}

A complete proof of this theorem can be found, for instance, in Chapter~7 of the recent book~\cite{BG06} by Bombieri and Gubler (who use a slightly different definition of height).

\section{Complexity of Algebraic Numbers}
\label{scomp}
Quite recently Adamczewski and Bugeaud applied the Subspace Theorem to  the long-standing problem of complexity of algebraic numbers. In particular, they proved transcendence of irrational automatic numbers. This will be the first topic of this talk.

We need some definitions. Let~$\AA$ be a finite set. We call it an \textit{alphabet}, and its elements will be referred to as \textit{letters}.   Let ${U=(u_1,u_2,u_3, \ldots)}$ be an infinite sequence of letters from~$\AA$. For every positive integer~$n$ we let ${\rho(n)=\rho_U(n)}$ the number of distinct $n$-words occurring as~$n$ successive elements of~$U$:
$$
\rho(n) = \bigl|\{u_ku_{k+1}\ldots u_{k+n-1}\ |\ k=1,2,3,\ldots\}\bigr|.
$$
Obviously, ${1\le \rho(n)\le |\AA|^n}$. The function~$\rho(n)$, defined on the set of natural numbers, is called the \emph{complexity function}, or simply \emph{complexity} of the sequence~$U$. 

Now let ${\alpha\in (0,1)}$ be a real number. For every integer ${b\ge 2}$ we can write the $b$-ary digital  expansion of~$\alpha$:
\begin{equation}
\label{ebary}
\alpha = u_1b^{-1}+u_2b^{-2}+ u_3b^{-3}+ \ldots, 
\end{equation}
where ${u_1,u_2,u_3, \ldots \in \{0,1, \ldots, b-1\}}$. One may ask about the complexity of the digital sequence ${(u_1,u_2,u_3, \ldots)}$. For instance, if~$\alpha$ is rational, then the expansion is (eventually) periodic, and the complexity function is bounded. Adamczewski and Bugeaud proved that the complexity function of the $b$-ary expansion of an irrational algebraic number is strictly non-linear.

\begin{rtheorem}{Adamczewski, Bugeaud}
\label{tcomp}
Let ${\alpha\in (0,1)}$ be an irrational algebraic number, and let ${b\ge 2}$ be an integer. Then the complexity function $\rho(n)$ of the $b$-ary expansion of~$\alpha$ satisfies
${\displaystyle\lim_{n\to \infty}\rho(n)/n=\infty}$. 
\end{rtheorem}

Previously, it was only known that ${\rho(n)-n\to +\infty}$, which  follows from the results of~\cite{FM97}. 

It is widely believed since the work of Borel~\cite{Bo09,Bo50}  that irrational algebraic numbers are \textit{normal}; that is, every $n$-word occurs in the $b$-ary expansion with the correct frequency~$b^{-n}$. In particular, one should expect that ${\rho(n)=b^n}$. This conjecture (let alone Borel normality) is far beyond the capabilities of the modern mathematics.

An important consequence of this theorem is transcendence of irrational \emph{automatic numbers}.  Recall that a finite automaton consists of the following elements:
\begin{itemize}
\item
the \emph{input alphabet}, which is usually the set of ${k\ge 2}$ digits ${\{0,1, \ldots, k-1\}}$;
\item
the set of \emph{states}~$\QQ$, usually a finite set of~$2$ or more elements,  with one element (called the \emph{initial state}) singled out;

\item the \emph{transition map} ${\QQ\times \{0,1, \ldots, k-1\}\to \QQ}$, which associates to every state a new state depending on the current input;

\item
the \emph{output alphabet}~$\AA$, together with the \emph{output map} ${\QQ\to \AA}$. 
\end{itemize}

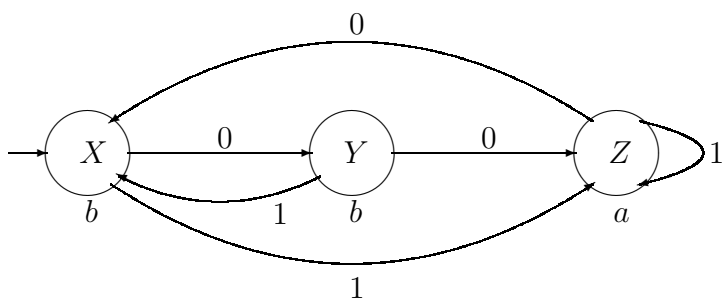
\begin{figure}
\begin{center}
\begin{picture}(260,100)(-30,-50)
\put(0,0){\circle{30}}
\put(100,0){\circle{30}}
\put(200,0){\circle{30}}
\put(15,0){\vector(1,0){70}}
\put(115,0){\vector(1,0){70}}
\qbezier(9,12)(100,72)(191,12)
\qbezier(9,-12)(100,-72)(191,-12)
\put(182,-18){\vector(3,2){10}}
\put(18,18){\vector(-3,-2){10}}
\qbezier(12,-9)(50,-28)(88,-9)
\put(21,-13){\vector(-2,1){10}}
\qbezier(209,12)(257,0)(209,-12)
\put(221,-9){\vector(-4,-1){13}}
\put(-3,-4){$X$}
\put(97,-4){$Y$}
\put(197,-4){$Z$}
\put(-1,-26){$b$}
\put(99,-26){$b$}
\put(199,-26){$a$}
\put(49,2){$0$}
\put(149,2){$0$}
\put(99,45){$0$}
\put(99,-55){$1$}
\put(70,-27){$1$}
\put(235,-4){$1$}
\put(-30,0){\vector(1,0){15}}
\end{picture}
\end{center}
\caption{A finite automaton with~$3$ states}
\label{faut}
\end{figure}

On Figure~\ref{faut} one can see an example of a finite automaton with inputs $0,1$, states $X,Y,Z$ with~$X$  the initial state, and outputs $a,b$. The transition map is given by the arrows, and the output map is ${X\mapsto b}$, ${Y\mapsto b}$ and ${Z\mapsto a}$.

An input stream for a finite automaton is a word in the input alphabet. Let us take the word $00100$. We start at the initial state~$X$ and the first input~$0$ moves us to the state~$Y$. The next  input~$0$ moves us further to~$Z$, and the third input~$1$ tells us to stay in~$Z$. With the fourth input~$0$ we return to~$X$, and with the final fifth input we end up in~$Y$. The output of~$Y$ is~$b$. Thus, the word $00100$ produces output~$b$. 

If we input consecutively  the binary expansion of natural numbers ${0,1,2,3, \ldots}$ written from right to left (that is, ${0,1, 01, 11, 001,\ldots}$), we obtain the sequence of outputs ${b,a,b,a,a,\ldots}$ called the \textit{automatic sequence} generated by the automaton from Figure~\ref{faut}. 

More generally, given an automaton with~$K$ inputs ${0,1,\ldots, k-1}$, the sequence generated by this automaton is the result of consecutive  inputs of $k$-ary expansions of natural numbers written from right to left. 

Probably, the most famous non-periodic automatic sequence is the \emph{Thue-Morse sequence}  ${0,1,1,0,1,0,0,1,\ldots}$, the $n$-th term being the parity of the sum of digits of the binary expansion of~$n$;  it is  generated by a finite automaton with~$2$ inputs,~$2$ states and~$2$ outputs. 

A real number ${\alpha \in (0,1)}$ is called \emph{automatic} if the digits of its $b$-ary expansion (for some ${b\ge 2}$) form an automatic sequence. 

For more information on automatic sequence see the book of Allouche  and Shallit~\cite{AS03}.

It is well-known (see, for instance,~\cite{Co72} or \cite[Section~10.3]{AS03}) that the complexity of an automatic sequence satisfies ${\rho(n)=O(n)}$. Hence Theorem~\ref{tcomp} implies the following remarkable result.

\begin{corollary}
An irrational automatic number is transcendental.
\end{corollary}

Probably, the first one to conjecture this was Cobham~\cite{Co68}. Sometimes this is referred to as the \textit{problem of Loxton and van der Poorten}, who obtained~\cite{LP82,LP88} several results in favor of this conjecture.

Adamczewski and Bugeaud deduce Theorem~\ref{tcomp} from a new transcendence criterion they obtained jointly with F.~Luca. The proof of this criterion relies on the Subspace Theorem. We say that the infinite sequence~$(u_n)$ has \emph{long repetitions} if there exist a real ${\eps>0}$, and infinitely many natural~$N$ such that  the word ${u_1u_2\ldots u_N}$ has two disjoint equal subwords of length exceeding~$\eps N$.

In symbols, the phrase ``the word ${u_1u_2\ldots u_N}$ has two disjoint equal subwords of length~$\ell$'' means the following: there exist~$k$ and~$n$ such that ${k+\ell\le n\le N+1-\ell}$ and 
$$
u_k=u_n, \quad u_{k+1}=u_{n+1},\quad \ldots,\quad u_{k+\ell-1}=u_{n+\ell-1}.
$$

\begin{rtheorem}{Adamczewski, Bugeaud, Luca}
\label{ttran}
Assume that for some ${b\ge 2}$ the $b$-ary expansion of  ${\alpha\in (0,1)}$ has long repetitions. Then~$\alpha$ is either rational or transcendental. 
\end{rtheorem}

In the introduction we remarked that the decimal expansion of an irrational algebraic number cannot have too long blocks of zeros (or too long periodic blocks), which is a relatively easy  consequence of the theorem of Ridout. Theorem~\ref{ttran} is a far-going   generalization of this observation.

Theorem~\ref{tcomp} is a consequence of Theorem~\ref{ttran}, due to the following simple lemma. 

\begin{lemma}
Assume that the complexity function of an infinite sequence~$(u_n)$ satisfies ${\displaystyle \liminf_{n\to\infty}\rho(n)/n<\infty}$. Then~$(u_n)$ has long repetitions.
\end{lemma}

\paragraph*{\textsc{Proof}} By the assumption, there exists ${\kappa>0}$ such that ${\rho(n)<\kappa n}$ for infinitely many~$n$. Fix such~$n$ and put ${N=\lceil(\kappa+1) n\rceil}$. By the box principle, the word ${u_1u_2\ldots u_N}$ contains two equal subwords of length~$n$. If they are disjoint, then we are done, because ${n\ge N/2(\kappa+1)}$.  Now assume they are not. This means that ${u_1u_2\ldots u_N}$ contains a subword ${W=ABC}$, where the words~$A$,~$B$ and~$C$ are non-empty and where $AB$ and~$BC$ are equal words of length~$n$.  

Since the words~$AB$ and~$BC$ are equal, we have ${W=AAB}$, which means that~$AA$ is a prefix\footnote{A \emph{prefix} of the word ${v_1\ldots v_m}$ is any of the words ${v_1\ldots v_s}$ with ${s\le m}$.} of~$W$. If ${\ell(AA)\le n}$ (where we denote by $\ell(X)$ the length of the word~$X$) then~$AA$ is a prefix of $AB$, which means that $AAA$ is a prefix of~$W$. Continuing by induction, we see that~$W$ has a prefix $\underbrace{A\ldots A}_k$, where ${k=\lfloor n/\ell(A)\rfloor+1}$ (in particular, ${k\ge 2}$ and ${k\ell(A)>n}$).  This implies that there are two disjoint words equal to $\underbrace{A\ldots A}_{\lfloor k/2\rfloor}$. Since ${k\ge 2}$ we have ${\lfloor k/2\rfloor\ge k/3}$, which implies that the length of these words is at least $n/3$. Hence the lemma is proved with ${\eps =1/6(\kappa+1)}$. \qed

\medskip

\paragraph*{\textsc{Proof of Theorem~\ref{ttran}}}
We assume that~$\alpha$ is algebraic and show that it is rational.
Write the $b$-ary expansion of~$\alpha$ as in~(\ref{ebary}).  By the hypothesis, there exist ${\eps>0}$ and infinitely many natural~$N$ such that the initial $N$-segment ${W_N=u_1\ldots u_N}$ has two disjoint subwords of length at least~$\eps N$. Fix one such~$N$. Then~$W_N$ has a prefix $ABCB$, where ${\ell(B)\ge \eps N}$ (the words~$A$ and~$C$ may be empty). Let~$\xi$ be the rational number with the eventually periodic $b$-ary expansion $ABCBCBC\ldots$. A straightforward calculation shows that 
$$
\xi=\frac M{b^r(b^s-1)},
$$
with ${M\in \Z}$, where  ${r=\ell(A)}$ is the length of the non-periodic part, and  ${s=\ell(BC)}$ is the length of the period. Notice that ${s+r=\ell(ABC)\le N}$ and that ${s\ge \ell(B)\ge \eps N}$. 

The main point of the proof is that~$\xi$ is a good rational approximation for~$\alpha$. Indeed, the first $\ell(ABCB)$ digits of the  $b$-ary expansions of~$\alpha$ and~$\xi$ coincide. Since ${\ell(ABCB)=r+s+\ell(B)\ge r+s+\eps N}$, we obtain
\begin{equation}
\label{eapp}
|\alpha-\xi|\le b^{-r-s-\eps N}
\end{equation}
This is not sufficient to get a contradiction with Roth's or Ridout's theorems, but, as we shall see, the Subspace Theorem will do the job.

Rewrite~(\ref{eapp}) as
\begin{equation}
\label{elfab}
\left|b^{r+s}\alpha-b^r\alpha-M\right|\le b^{-\eps N}.
\end{equation}
Now it is the time to define the data for the Subspace Theorem: the set~$S$ of prime numbers and  the linear forms $L_{i,p}$. 
Let~$S$ consist of the infinite prime and all the prime divisors of~$b$. Further, for ${p\in S}$ we define the linear forms~$L_{1,p}$,~$L_{2,p}$ and~$L_{3,p}$ in variables ${\bx=(x_1,x_2,x_3)}$ as follows. For ${p=\infty}$ we put 
$$
L_{1,\infty}(\bx)=x_1, \qquad L_{2,\infty}(\bx)=x_2, \qquad L_{3,\infty}(\bx)=\alpha x_1 -\alpha x_2 -x_3.
$$
And for ${p < \infty}$ we put ${L_{i,p}(\bx)=x_i}$ for ${i=1,2,3}$. 

We put ${\bx= \left(b^{r+s}, b^s, M\right)}$. Since ${\xi\in (0,1)}$,  we have ${|M|\le b^{r+s}}$. Thus,
\begin{equation}
\label{enorm}
\|\bx\|\le b^{r+s}\le b^N.
\end{equation}
Now we have
$$
\prod_{p\in S}\prod_{i=1}^3\left|L_{i,p}(\bx)\right|_p= \prod_{p\in S}|b^r|_p \prod_{p\in S}|b^{r+s}|_p
\prod_{\genfrac{}{}{0pt}{}{p\in S}{p\ne\infty}}|M|_p \left|b^{r+s}\alpha-b^r\alpha-M\right|_\infty
$$
By our definition of~$S$ and the product formula we have ${\prod_{p\in S}|b|_p=1}$. Further, since ${M\in \Z}$, we have 
${|M|_p\le 1}$ for each ${p\ne \infty}$. It follows that 
\begin{equation}
\label{eiab}
\prod_{p\in S}\prod_{i=1}^3\left|L_{i,p}(\bx)\right|_p\le\left|b^{r+s}\alpha-b^r\alpha-M\right|_\infty \le b^{-\eps N} \le \|\bx\|^{-\eps}.
\end{equation}
(We used~(\ref{elfab}) and~(\ref{enorm}).)

We can repeat this argument for infinitely many~$N$ and find vectors ${\bx=\bx(N)}$ satisfying~(\ref{eiab}). Moreover, recall that  ${s(N)\ge  \eps N}$, whence ${s(N)\to \infty}$ as ${N\to \infty}$, which means that among the vectors ${\bx=\bx(N)}$ infinitely many are distinct. Theorem~\ref{tschl} implies that these vectors $\bx(N)$ lie on finitely many planes of the space~$\Q^3$. Hence infinitely many of them lie on the same plane; that is, there exist ${\lambda, \mu, \nu \in \Q}$, not all~$0$ such that for infinitely many~$N$ we have
\begin{equation}
\label{elin}
\lambda b^{r(N)} +\mu b^{r(N)+s(N)} +\nu M(N)=0.
\end{equation}
Moreover, ${\nu\ne 0}$ because ${s(N)\to \infty}$. Dividing~(\ref{elin}) by ${b^{r(N)}\left(b^{s(N)}-1\right)}$, we obtain 
$$
\frac\lambda{b^{s(N)}-1} + \mu \frac{b^{s(N)}}{b^{s(N)}-1} +\nu \xi(N)=0.
$$
Sending~$N$ to infinity, we conclude that ${\mu+\nu\alpha =0}$, whence ${\alpha \in \Q}$. The theorem is proved. \qed

\medskip

As the reader could have noticed, it is quite irrelevant for the proof that the ``digits'' ${u_1,u_2, \ldots}$ belong to the set ${\{0,1, \ldots, b-1\}}$. In fact, any finite set of rational, or even algebraic numbers would do. Also,~$b$ is not obliged to be a rational integer; one can assume it to be any Pisot or Salem number\footnote{\label{fpis} A real algebraic number ${\beta>1}$ is called \emph{Pisot number} if all its conjugates (except~$\beta$ itself) lie inside the unit disk of the complex plane; it is called \emph{Salem number} if they lie inside or on the boundary of the unit disk.}. Thus, the result of Adamczewski and Bugeaud   in the most general form sounds as follows: let ${u_1,u_2, \ldots}$ be a sequence of algebraic numbers with finitely many distinct terms, and with long repetitions, and let~$\beta$ be a Pisot or Salem number; then either 
${\alpha= u_1\beta^{-1}+ u_2\beta^{-2}+\ldots}$ belongs to the number field generated by~$\beta$ and by the ``digits'' ${u_1,u_2, \ldots}$, or~$\alpha$ is transcendental.

In~\cite{AB05,AB07a} Adamczewski and Bugeaud exploit a different notion of complexity, based on continued fractions rather than $b$-ary expansions, and obtain several results in the same spirit. 

The reader may consult Waldschmidt's survey~\cite{Wa07} for more information on the Diophantine analysis of symbolic sequences. 

Remark in conclusion that Adamczewski and Bugeaud were not the first to apply the Subspace Theorem in the transcendence; in \cite{Mi77,BP87,TZ99,CZ02a} it was used to prove transcendence of certain infinite sums.  The argument of Troi and Zannier~\cite{TZ99} is quite similar to that of Adamczewski and Bugeaud. See also \cite{CZ06a} for a more recent application.

\begin{sloppypar}
\section{Diophantine Equations with Power Sums}
\label{spowsum}
In 1984 M.~Laurent~\cite{La84} applied the Subspace Theorem to study the solutions ${\bx =(x_1, \ldots, x_r)\in \Z^r}$ of a polynomial-exponential equation 
	\begin{equation}
	\label{epolexp}
	\sum_{i=1}^NP_i(\bx)\ba_i^\bx=0,
	\end{equation}
	where ${P_1,\ldots,P_N \in \barQ[\bx]}$,  ${\ba_1, \ldots, \ba_N\in \barQ^{r}}$ and ${\ba^\bx:=a_1^{x_1}\cdots a_{r\vphantom1}^{x_r}}$. (Using a specialization argument, one can replace~$\barQ$ by any field of characteristic~$0$.) His results imply, in particular, that, under certain natural condition, the following holds: with finitely many exceptions, every solution of~(\ref{epolexp}) is also a solution of a ``strictly shorter'' equation  ${\sum_{i\in I}P_i(\bx)\ba_i^\bx=0}$, where~$I$ is a proper subset of ${\{1, \ldots, N\}}$. The above mentioned condition is the following:
	the only ${\bx\in \Z^r}$ satisfying  ${\ba_i^\bx=\ba_j^\bx}$ for all $i,j$ is ${\bx =(0,\ldots, 0)}$. 
\end{sloppypar}

While the theorem of Laurent does not (and cannot) imply ultimate finiteness in general, it allows one to establish it in many special cases, usually by induction in~$N$ and/or elimination. 

However, there are interesting polynomial-exponential equations for which the theorem of Laurent does not yield anything non-trivial. One of the simplest is ${a^n+b^n=P(x)}$ in ${x,n\in \Z}$, where~$P$ is a polynomial. (The equation ${a^n=P(x)}$ can be analyzed, for instance, by Baker's method.) For this equation ${r=2}$ and the vectors~$\ba_i$ are  $(a,1)$, $(b,1)$ and $(1,1)$. For any such~$\ba_i$ and for any ${\bx=(0,x)}$ we have ${\ba_i^\bx=(1,1)}$, so that Laurent's  condition is not satisfied.

Corvaja and Zannier studied these and more general equations in the important and largely underestimated article~\cite{CZ98}, as well as in the later article~\cite{CZ02a}. Let us introduce some terminology. Call \emph{power sum} an expression of the form \begin{equation}
\label{epowsum}
u(n)=b_1a_1^n+\cdots+b_ma_m^n, 
\end{equation}
where ${a_1, \ldots, a_m}$ (the  \emph{roots}) and ${b_1, \ldots, b_m}$ (the \emph{coefficients}) are complex numbers. Power sums can be viewed as a particular case of \emph{linear recurrence sequences} 
$$
u(n)=b_1(n)a_1^n+\cdots+b_m(n)a_m^n,
$$ 
where ${b_1(n), \ldots , b_m(n)}$ are polynomials in~$n$; one may say that  \emph{power sums are linear recurrences with simple roots}.

If the roots and the coefficients  belong to a ring~$A$, then we call~(\ref{epowsum}) an \emph{$A$-power sum}, or a \emph{power sum over~$A$}.

Let ${P(x,y)\in \Q[x,y]}$ be an irreducible polynomial with ${\deg_yP\ge 2}$. Corvaja and Zannier studied the equation ${P(u(n), y)=0}$, where~$u$ is a power sum. They were motivated by a question of Yasumoto about \emph{universal Hilbert sets}, that is, sets~$A$ of rational integers with the following property: 
\begin{center}
\begin{tabular}{rl}
(UHS)&\parbox{\tmtwo}{for any irreducible (over~$\Q$) polynomial  ${P(x,y)\in \Q[x,y]}$,  the specialized polynomial ${P(a,y)\in \Q[y]}$ is irreducible for all but finitely many ${a\in A}$.}
\end{tabular}
\end{center}
Informally, a universal Hilbert  set proves the Hilbert irreducibility theorem for every polynomial, and with finitely many exceptions. 

A well-known elementary Galois-theoretic argument (see, for instance, \cite[Section~2]{Bi96}) implies that~$A$ is a universal Hilbert set if and only if it has the following formally weaker property: 
\begin{center}
\begin{tabular}{rl}
(UHS${}'$)&\parbox{\tmtwo}{for any absolutely irreducible ${P(x,y)\in \Q[x,y]}$ with ${\deg_yP\ge 2}$ the equation ${P(a, y)=0}$ has only finitely many solutions $(a,y)$ with ${a\in A}$ and ${y\in \Q}$.}
\end{tabular}
\end{center}

Existence of universal Hilbert sets was shown by Gilmore and Robinson~\cite{GR55}, and the first  explicit example was suggested by Sprindzhuk~\cite{Sp81} (see \cite{Bi96,DZ98,Ya87,Za96} for further examples). Yasumoto~\cite{Ya87} asked whether ${\{2^n+3^n\}}$ is a universal Hilbert set. Dèbes and Zannier~\cite{DZ98} managed to prove, using the theorem of Ridout, that ${\{2^n+5^n\}}$ is a universal Hilbert set, but their argument fails for ${\{2^n+3^n\}}$. In~\cite{CZ98} this problem is solved, and even a much stronger result is obtained: values of any power sum ${b_1a_1^n+\cdots+b_ma_m^n}$ with multiplicatively independent
${a_1, \ldots, a_m}$ form a universal Hilbert set (with ${m\ge 2}$ and ${b_1, \ldots, b_m\ne 0}$).

Another motivation for~\cite{CZ98} was the celebrated problem of Pisot. A power series ${f(t)=\sum_{n=0}^\infty u(n)t^n}$ is called the \emph{Hadamard $q$-th power} of the series ${g(t)=\sum_{n=0}^\infty v(n)t^n}$ if ${u(n)=v(n)^q}$ for ${n=0,1,\ldots}$; in this case the latter series is called an \emph{Hadamard $q$-th root} of the former.

Let $f(t)$ be a rational power series (that is, a power series expansion of a rational function in~$t$) with coefficients in~$\Q$, and let~$q$ be a positive integer. Assume that $f(t)$ is the Hadamard $q$-th power of another series with coefficients in~$\Q$. Pisot conjectured that in this case $f(t)$ is the Hadamard $q$-th power of another \emph{rational} power series (with coefficients in~$\Q$). 

Since ${f(t)=\sum_{n=0}^\infty u(n)t^n}$ is a rational power series if and only if the coefficients~$u(n)$ form a linear recurrence sequence, Pisot's conjecture can be stated as follows: assume that $\{u(n)\}$ is a linear recurrence sequence of rational numbers, such that every~$u(n)$ is a $q$-th power in~$\Q$; then ${u(n)=v(n)^q}$ for all~$n$, where $v(n)$ is another linear recurrence sequence of rational numbers. 

Zannier~\cite{Za00} proved Pisot's conjecture by a method independent of the Subspace Theorem. Now, let us ask a more difficult question: assume that
\begin{equation}
\label{erefpis}
\text{$u(n)$ is a $q$-th power in~$\Q$ for \emph{infinitely many}~$n$;} 
\end{equation}
what can one say about the linear recurrence~$u$? Since the work of Corvaja and Zannier applies to the particular equation ${u(n)-y^q=0}$, it answers this question in the special case when~$u$ is a power sum (over~$\Q$).   It turns out that, while~$u$ itself is not obliged to be a $q$-th power of another $\Q$-power sum, this is true for the power sum obtained from~$u$ by letting~$n$ run through an arithmetical progression (see Corollary~\ref{cps}).

Below, we give a complete proof of this particular case of the theorem of Corvaja and Zannier. We shall also state the general theorem and sketch its proof.

\subsection{Refined Pisot's Conjecture for Power Sums}

The main result of Corvaja and Zannier concerns $\Q$-power sums with \emph{positive} roots. For these power sums~(\ref{erefpis}) implies that~$u$ is a $q$-th power of another power sum, but over~$\bar \Q$. More precisely, we have the following.
 
\begin{rtheorem}{Corvaja, Zannier}
\label{tpos}
Let~$u$ be a $\Q$-power sum with positive roots, and let~$q$ be a positive integer. Assume that $u(n)$ is a $q$-th power for infinitely many ${n\in \Z}$. Then ${u(n)=a^{n+r}v(n)^q}$ for all ${n\in \Z}$, where~$a$ is a non-zero rational number,~$r$ is an integer and~$v$ is a $\Q$-power sum. In particular,~$u$ is a $q$-th power of  $\bar \Q$-power sum. 
\end{rtheorem}

\begin{corollary}
\label{cps}
Let~$u$ be a $\Q$-power sum, and let~$q$ be a positive integer. Assume that $u(n)$ is a $q$-th power for infinitely many ${n\in \Z}$. Then there exist positive integers~$Q$ and~$R$ and a $\Q$-power sum~$w$ such that ${u(Qn+R)=w(n)^q}$ for all ${n\in \Z}$.
\end{corollary}

In other words, though~$u$ itself is not necessarily a $q$-th power of a $\Q$-power sum, the power sum obtained from~$u$ by letting~$n$ run through a certain arithmetical progression is. 

If~$u$ has positive roots then the corollary is immediate, with ${Q=q}$. In the general case one should consider the power sums $u(2n)$ and ${u(2n+1)}$, both having positive roots, and the corollary follows with ${Q=2q}$.  

\paragraph*{\textsc{Proof of Theorem~\ref{tpos}}}
We may assume that~$u(n)$ is a $q$-th power  for infinitely many 
\emph{positive} integers~$n$, replacing $u(n)$ by ${u(-n)}$, if necessary. 

Write ${u(n)=b_0a_0^n+\cdots+b_ma_m^n}$, where the roots ${a_0, \ldots, a_m}$ are positive rational numbers written in the decreasing order,  so that ${a_0>a_1>\ldots>a_m>0}$.  

Assume first  that ${a_0=1}$. Putting ${b=b_0}$ and ${c_k=b_k/b}$, we  write  
$$
u(n)=b(1+z(n)) 
$$
with ${z(n)=c_1a_1^n+\cdots +c_ma_m^n}$. Since the roots of the power sum~$z$  are strictly smaller than~$1$, we have\footnote{In this proof ``$\ll$'', ``$\gg$'' and $O(\cdot)$ imply  constants depending on the power sum~$u$ and on the parameter~$\Lambda$ defined below, but independent of~$n$.}, ${|z(n)|\ll\theta^n}$ for some ${\theta \in (0,1)}$. Since $u(n)$ is infinitely often a rational $q$-th power, we have ${b>0}$ when~$q$ is even. We may assume that ${b>0}$ when~$q$ is odd as well,  replacing~$u$ by~${-u}$, if necessary. Thus, for big positive~$n$ we have ${u(n)>0}$, which implies that~$u(n)$ has exactly one positive $q$-th root; we denote it by $y(n)$. For sufficiently large~$n$ we can express $y(n)$ using the binomial power series:
\begin{equation}
\label{eroot}
y(n)= b^{1/q}\sum_{\ell=0}^{\Lambda-1}\binom{1/q}\ell z(n)^\ell+O(\theta^{n\Lambda}),
\end{equation}
where the parameter~$\Lambda$  will be specified later. 

The sum in~(\ref{eroot}) can be expressed as ${\beta_1\alpha_1^n+ \cdots+\beta_\mu\alpha_\mu^n}$, where ${\alpha_1, \ldots, \alpha_\mu}$ are pairwise distinct. Since ${\alpha_1, \ldots, \alpha_\mu}$ are multiplicative combinations of ${a_1, \ldots, a_m}$,  they are positive rational numbers. Thus, we have
$$
\left|y(n)-b^{1/q}\sum_{k=1}^\mu\beta_k\alpha_k^n\right|\ll\theta^{n\Lambda}.
$$
Now we are in a position to apply the Subspace Theorem. We let~$S$ to be a finite set of prime numbers, including the infinite prime, such that the numbers ${a_1, \ldots, a_m}$ are $S$-units, and ${b_0, \ldots, b_m}$ are $S$-integers. Then $u(n)$ is an $S$-integer for every~$n$, and so is ${y(n)=u(n)^{1/q}}$ (as soon as ${y(n)\in \Q}$). Also, the numbers ${\alpha_1, \ldots, \alpha_\mu}$ are $S$-units, being multiplicative combinations of ${a_1, \ldots, a_m}$. 

Next, for every ${p\in S}$ we define ${\mu+1}$ independent linear forms  in ${\mu+1}$ variables as follows.  For ${p=\infty}$ we put 
$$
L_{0,\infty}(\bx)= x_0-b^{1/q}\sum_{k=1}^\mu\beta_kx_k, \qquad
L_{k, \infty}(\bx) = x_k \quad (k=1, \ldots, \mu).
$$
And for a finite ${p\in S}$ we put ${L_k(\bx)=x_k}$ for ${k=0, \ldots, \mu}$. 

Now let~$n$ be such that ${y(n)\in \Q}$. Then  ${\bx=\bx(n)=\left(y(n), \alpha_1^n, \ldots, \alpha_\mu ^n\right)}$ is a vector with $S$-integer coordinates. We have
\begin{equation}
\label{epr1}
\prod_{p\in S}\prod_{k=0}^r\left|L_{k,p}(\bx)\right|_p=\left|y(n)-b^{1/q}\sum_{k=1}^\mu\beta_k\alpha_k^n\right|  \prod_{\genfrac{}{}{0pt}{}{p\in S}{p\ne \infty}}|y(n)|_p  \prod_{k=1}^r\prod_{p\in S}\left|\alpha_k^n\right|_p\ll\theta^{\Lambda n}H(y(n)).
\end{equation}
Indeed, the product formula implies that
${\prod_{p\in S}|\alpha_k|_p=1}$ (because the numbers~$\alpha_k$ are $S$-units), which means that the double product is~$1$. Also, the first product is bounded by $H(y(n))$, because $y(n)$ is an $S$-integer.

An obvious calculation shows that the height of the rational number $u(n)$ is $e^{O(n)}$. Since ${y(n)^q=u(n)}$, we have ${H(y(n))=H(u(n))^{1/q}=e^{O(n)}}$. It follows that the right-hand side of~(\ref{epr1}) is bounded by ${C^n\theta^{\Lambda n}}$, where the constant~$C$ depends only on the power sum~$u$. 

Now specify the parameter~$\Lambda$ to have ${C\theta^\Lambda\le 1/2}$. We obtain
$$
\prod_{p\in S}\prod_{k=0}^r\left|L_{k,p}(\bx)\right|_p\ll 2^{-n}.
$$
A routine estimate gives ${H(\bx)\le e^{O(n)}}$. We finally obtain
\begin{equation}
\label{eschl}
\prod_{p\in S}\prod_{k=0}^r\left|L_{k,p}(\bx)\right|_p< H(\bx)^{-\eps}
\end{equation}
with some ${\eps>0}$ (depending only on~$u$). 

By the assumption, there exist infinitely many positive integers~$n$ such that ${y(n)\in \Q}$. Hence~(\ref{eschl}) has infinitely many solutions in $S$-integer vectors ${\bx=\bx(n)}$. By Theorem~\ref{tschl}${}'$, all these solutions belong to finitely many proper subspaces of~$\Q^{\mu+1}$. It follows that infinitely many vectors $\bx(n)$ belong to the same proper subspace. In other words, there exist rational numbers ${\gamma_0,\ldots,  \gamma_\mu }$, not all~$0$, such that 
$$
\gamma_0y(n) = b^{1/q}\left(\gamma_1\alpha_1^n+ \cdots + \gamma_\mu \alpha_\mu ^n\right).
$$
If ${\gamma_0=0}$ then ${\gamma_1\alpha_1^n+ \cdots + \gamma_\mu \alpha_\mu ^n}$ would vanish for infinitely many~$n$, which is impossible because ${\alpha_1, \ldots, \alpha_\mu }$ are pairwise distinct positive numbers. Thus, ${\gamma_0\ne 0}$, and we may assume that ${\gamma_0=1}$.

We have shown that for infinitely many~$n$ we have ${y(n)\in \Q}$ and 
$$
y(n)=  b^{1/q}\left(\gamma_1\alpha_1^n+ \cdots + \gamma_\mu \alpha_\mu ^n\right).
$$
Since ${\gamma_1\alpha_1^n+ \cdots + \gamma_\mu \alpha_\mu ^n\ne 0}$ for large~$n$, we have   
$$
b^{1/q}= \frac{y(n)}{\gamma_1\alpha_1^n+ \cdots + \gamma_\mu \alpha_\mu ^n} \in \Q.
$$
Thus, for infinitely many~$n$ we have ${y(n)=v(n)}$, where~$v$ is a $\Q$-power sum with positive roots. Since ${u(n)-v(n)^q}$ is a power sum with positive roots as well, it can vanish infinitely often only if it vanishes identically. Thus, ${u(n)=v(n)^q}$. This proves the theorem in the special case ${a_0=1}$. 

The general case easily reduces to the special one. For some~$r$ there exist infinitely many positive integers~$n$, congruent to~${-r}$ modulo~$q$ such that $u(n)$ is a $q$-th power in~$\Q$. Replacing $u(n)$ by ${a_0^{-n-r}u(n)}$, we reduce the general case to the case ${a_0=1}$, already treated. \qed

\subsection{The General Equation}

And here is the general theorem of Corvaja and Zannier.

\begin{rtheorem}{Corvaja, Zannier~\cite{CZ02a}}
\label{tgencz}
Let~$u$ be a $\Q$-power sum with positive roots, let~$S$ be a finite set of primes including the infinite prime, and  let ${P(x,y)\in \Q[x,y]}$ be a polynomial non-constant in~$y$. Assume that the equation ${P(u(n),y)=0}$ has infinitely many solutions in integers~$n$ and $S$-integers~$y$.  Then there exists a $\bar \Q$-power sum~$v$ with positive real coefficients such that ${P(u(n), v(n))=0}$ for all ${n\in \Z}$. 
\end{rtheorem}

\paragraph*{\textsc{Proof} (a sketch)}
As above, we may assume that there are infinitely many solutions with positive~$n$. When  ${n\to +\infty}$ we have ${u(n)\to b\in \Q\cup \{-\infty,+\infty\}}$. Replacing $u(n)$ by $-u(n)$ we may exclude the ${-\infty}$, and upon replacing ${u(n)}$ by ${u(n)-b}$ we may assume that in the finite case the limit is~$0$. Thus, ${\displaystyle\lim_{n\to+\infty}u(n)\in \{0, +\infty\}}$. 

Assume that ${\displaystyle\lim_{n\to+\infty}u(n)=0}$. Then\footnote{Here implicit  constants may depend on the power sum~$u$, the polynomial $P(x,y)$ and the parameter~$\Lambda$ defined below, but not on~$n$.} ${u(n)\ll \theta^n}$ with some ${\theta\in (0,1)}$.
By the assumption, for infinitely many positive integers~$n$ there exists an $S$-integer $y(n)$ such that 
${P(u(n), y(n))=0}$. Let
\begin{equation}
\label{epui}
Y_i(x)=\sum_{k=-\kappa_i}^\infty c_{ki}x^{k/e_i} \qquad (i=1, \ldots, \deg_yP)
\end{equation}
be the Puiseux expansion of the algebraic function~$y$ at~$0$, the coefficients~$c_{ki}$ being algebraic numbers. Since ${u(n)\to 0}$, for large~$n$ all the series~(\ref{epui})  converge at ${x=u(n)}$, and one of the sums $Y_i(u(n))$ is $y(n)$. We fix~$i$ for which ${Y_i(u(n))=y(n)}$ infinitely often, and omit the index~$i$ in the sequel. Thus, for infinitely many positive integers~$n$ we have
$$
y(n)=\sum_{k=-\kappa}^\infty c_ku(n)^{k/e}.
$$
Truncating the series, we find
$$
y(n)=\sum_{k=-\kappa}^{\Lambda e-1}c_ku(n)^{k/e} + O(\theta^{\Lambda n}).
$$

Now write ${u(n)=ba^n(1+z(n))}$, where~$a$ is the biggest root of~$u$. Redefining~$\theta$, we may assume that ${z(n)\ll\theta^n}$ for positive~$n$. Replacing each $u(n)^{k/e}$ by 
$$
b^{1/e}a^{n/e}\sum_{j=0}^{\Lambda-1}\binom{k/e}jz(n)^j +O(\theta^{\Lambda n}),
$$
we obtain 
${y(n)= \beta_1\alpha_1^n+\cdots+\beta_\mu \alpha_\mu^n + O(\theta^{\Lambda n})}$,
where ${\alpha_1, \ldots, \alpha_\mu}$ are positive real algebraic numbers, and ${\beta_1, \ldots, \beta_\mu}$ are algebraic numbers.

Now applying the Subspace Theorem in the same way as we did in the proof of Theorem~\ref{tpos}, we find that ${y(n)=v(n)}$ for infinitely many~$n$, where~$v$ is a $\bar \Q$-power sum with positive real roots.
Then ${P(u(n), v(n))}$ is a $\bar \Q$-power sum with positive real roots, which vanishes at infinitely many~$n$. Hence it vanishes identically. 

The case ${u(n)\to+\infty}$ is treated similarly, the Puiseux expansions at zero being replaced by those at infinity. \qed

\medskip

Among other consequences of this theorem, we have the following result mentioned above.

\begin{rcorollary}{Corvaja, Zannier}
\label{cuhs}
Let ${u(n)=b_1a_1^n+\cdots+b_ma_m^n}$ be a $\Q$-power sum. Assume that ${m\ge 2}$ and that the roots ${a_1, \ldots, a_m}$ are multiplicatively independent. Then ${\{u(n)\}}$ is a universal Hilbert set.
\end{rcorollary}

To prove the corollary, we need a purely algebraic lemma. Let~$K$ be a field of characteristic~$0$ and let~$\Gamma$ be a multiplicatively written torsion-free abelian group. Then the group ring ${K[\Gamma]}$ is an integral domain.

\begin{sloppypar}
	\begin{lemma}
	\label{lgam}
	 In the ring ${K[\Gamma]}$ consider an element 
	${u=b_1\gamma_1+\cdots +b_m\gamma_m}$, where ${b_1, \ldots, b_m\in K^\ast}$  and  ${\gamma_1, \ldots, \gamma_m}$ are multiplicatively independent elements of~$\Gamma$. Assume that ${m\ge 2}$. Then the ring $K[u]$ is integrally closed in $K[\Gamma]$. 
	\end{lemma}
\end{sloppypar}

Since this the lemma has nothing to do with our main subject, we prove it in the addendum to this section. 

\paragraph*{\textsc{Proof of Corollary~\ref{cuhs}}}
We apply the lemma with ${K=\C}$ and with~$\Gamma$ consisting of the functions ${\Z\to \R}$ defined by ${n\mapsto a^n}$ with a positive real~$a$. Then ${\C[\Gamma]}$ is exactly the ring of power sums with complex coefficients and positive real roots.

Now let ${u(n)=b_1a_1^n+\cdots+b_ma_m^n}$ be as in the corollary. We may assume that the roots ${a_1, \ldots, a_m}$ are positive, considering separately ${u(2n)}$ and ${u(2n+1)}$. By the lemma, the ring ${\C[u]}$ is integrally closed in $\C[\Gamma]$.

If $\{u(n)\}$ is not a universal Hilbert set then there exists a $\Q$-irreducible polynomial ${P(x,y)\in \Q[x,y]}$ with ${\deg_yP\ge 2}$ such that ${P(u(n), y)=0}$ has infinitely many solutions in ${n\in \Z}$ and ${y\in \Q}$. We may assume the polynomial~$P$ absolutely irreducible\footnote{It is well-known and easy to show that if $P(x,y)$ is $\Q$-irreducible but $\C$-reducible then the equation ${P(x,y)=0}$ can have only finitely many solutions in ${x,y\in \Q}$.} and monic\footnote{Replace ${P(x,y)=a_q(x)y^q+\cdots+a_1(x)y+1}$ by ${a_q(x)^{q-1}P\left(x, y/a_q(x)\right)}$.} in~$y$. Since~$P$ is monic, there exists a finite set of primes~$S$ such that for all solutions $(n, y)$ as above, the number~$y$ is an $S$-integer. Applying Theorem~\ref{tgencz}, we find a power sum~$v$ with positive real roots such that ${P(u(n),v(n))=0}$. Since the polynomial $P(x,y)$ is absolutely irreducible, $y$-monic  and of $y$-degree at least~$2$, the ring ${\C[u,v]}$ is a non-trivial integral extension of $\C[u]$. Hence ${\C[u]}$ is not integrally closed in $\C[\Gamma]$, a contradiction. \qed

\medskip

In fact, Corvaja and Zannier prove more. For instance, using Siegel's theorem (see Section~\ref{sipon}), they show\footnote{In~\cite{CZ98} they consider only power sums with integer roots, but the argument extends to rational roots without trouble.} the following: \textsl{in the set-up of Theorem~\ref{tgencz} assume that~$P$ is $\Q$-irreducible and ${\deg_yP\ge 2}$; then  either ${u=f(v)}$, where~$v$ is another power sum and~$f$ is a polynomial of degree at least~$2$, or the roots of~$u$ generate a cyclic multiplicative group (that is ${u(n)=b_1a^{\nu_1n}+\cdots+b_ma^{\nu_mn}}$ with some ${a\in \Q^\ast}$ and ${\nu_1, \ldots, \nu_m\in \Z}$).} This implies further examples of universal Hilbert power sums, like ${2^n+3^n+6^n}$, etc. 

To conclude,  we briefly discuss power sums over number fields.
Theorems~\ref{tpos} and~\ref{tgencz} stay true, with almost the same proof, if the assumption \textsl{the roots of~$u$ are positive} is replaced by \textsl{the roots of~$u$ generate a torsion-free multiplicative abelian group}. One may attempt to extend Theorems~\ref{tpos} and~\ref{tgencz}, with this more general assumption, to $K$-power sums, with an arbitrary number field~$K$. Unfortunately, this is done only under a certain technical assumption about our power sum. We say that a $K$-power sum~$u$  has an \emph{upper} (respectively, \emph{lower}) \emph{dominant root} if there exists a root~$a$ of~$u$ and an absolute value ${v\in M_K}$  such that ${|a|_v>|a'|_v}$ (respectively, ${|a|_v<|a'|_v}$) for any other root~$a'$. 

Now let~$u$ be a $K$-power sum satisfying the following two conditions: \textsl{the roots of~$u$ generate a torsion-free multiplicative abelian group}, and\textsl{~$u$ has both an upper dominant root and a lower dominant root.} Then~$u$ satisfies both the analogues of Theorems~\ref{tpos} and~\ref{tgencz} with~$\Q$ replaced by~$K$ (with very similar proofs). 

The existence of a ``dominant root'' is immediate\footnote{provided the roots generate a torsion-free abelian group} if the field~$K$ has at least one real embedding, but it may fail already for ${K=\Q(i)}$: the power sum 
$$
u(n) =(8+i)^n+(8-i)^n+(2+i)^n+(2-i)^n
$$
has no upper dominant root.

Suppressing  the ``dominant root'' assumption looks a difficult problem. It seems that at least one cardinal new idea is needed to handle power sums without dominant roots. See, however,~\cite{CZ02b}.

\subsection*{Addendum: Proof of Lemma~\ref{lgam}}

We may assume  that~$\Gamma$ is a division group; moreover, since it is torsion-free, every ${\gamma\in \Gamma}$ has a well-defined ``$n$-th root''~$\gamma^{1/n}$ for any non-zero integer~$n$.  It suffices to prove that $K[u]$ is integrally closed in the ring $K[\Delta]$, for any finitely generated subgroup~$\Delta$ of~$\Gamma$, containing ${\gamma_1,\ldots, \gamma_m}$. Replacing~$\Delta$ by a bigger finitely generated subgroup, we may assume that it has a free $\Z$-basis consisting of ${\gamma_1^{1/n}, \ldots, \gamma_m^{1/n}}$ (with some positive integer~$n$) and, perhaps,  several more elements of~$\Gamma$. 

We have reduced the lemma to the following statement.

\begin{proposition}
Let ${R=K[x_1, \ldots, x_r]}$ be the polynomial ring over a field~$K$ (of characteristic~$0$) and let~$n$ be a positive integer. Consider ${u=b_1x_1^n+ \cdots+b_mx_m^n\in R}$, where ${2\le m\le r}$ and ${b_1, \ldots, b_m\in K^\ast}$. Then $K[u]$ is integrally closed in~$R$. 
\end{proposition}

\paragraph*{\textsc{Proof}} 
We may assume that~$K$ is algebraically closed  and, by a linear change of variables we may assume that ${b_1= \ldots= b_m=1}$, so that ${u=x_1^n+ \cdots+x_m^n}$. Let~$\OO$ be the integral closure of $K[u]$ in~$R$.   We want to prove that ${\OO=K[u]}$.

The quotient field of~$\OO$ is contained in the purely transcendental field ${K(x_1, \ldots, x_r)}$. By the theorem of Luroth (see Remark~\ref{rluro}) it itself must be purely transcendental. Thus, we may write this quotient field  as  $K(v)$, and the generator~$v$ may be chosen in the ring~$\OO$. We have ${u=P(v)}$, where, a priori, $P(X)$ is a rational function over~$K$. Since both~$u$ and~$v$ are polynomials in ${x_1, \ldots, x_r}$,  the rational function $P(X)$ must be  a  polynomial.

Specializing ${x_1=t}$, ${x_2=\ldots=x_r=0}$, we obtain ${t^n=P(Q(t))}$, where~$Q(t)$ is a polynomial over~$K$. It follows that ${P(X)=aX^\nu}$ for some positive integer~$\nu$ and some ${a\in K^\ast}$. Specializing ${x_1=t}$, ${x_2=1}$, ${x_3=\ldots=x_r=0}$
(it is here where we use the assumption ${m\ge 2}$), we conclude that ${t^n+1}$ is a $\nu$-th power of yet another polynomial in~$t$, which is possible only if ${\nu =1}$. Thus, ${u=av}$, which proves the proposition.\qed

\begin{remark}
\label{rluro}
We use here a slightly non-traditional form of Luroth's theorem: if ${K\subset L\subset \Omega}$ is a tower of fields of characteristic~$0$, with~$K$ algebraically closed,~$\Omega$ purely transcendental over~$K$ and~$L$ of transcendence degree~$1$ over~$K$, then~$L$ is purely transcendental. In standard textbooks one usually assumes that~$\Omega$ is of transcendence degree~$1$ as well.

However, our ``more general'' version of Luroth's theorem easily follows from the traditional one. Indeed, geometrically, the ``traditional'' version means the following:  if an algebraic curve~$C$ admits a non-constant rational dominant map ${\P^1\to C}$, then it is isomorphic to~$\P^1$. And in our version~$\P^1$ should be replaced by~$\P^r$. But if a curve admits a non-constant dominant map from a projective space, then it also admits one from the projective line.
\end{remark}

\section{Integral Points}
\label{sipon}

\subsection{Integral Points on Curves}
It is well-known that a  binary Diophantine equation ${P(x,y)=0}$ of degree~$1$ or~$2$ has infinitely many solutions in integers unless it has an ``obvious'' reason (local obstruction) for having finitely many. Siegel proved~\cite{Si29}, relying on the already mentioned work of 
A.~Thue~\cite{Th09}, that an  equation  of degree~$3$ or higher must have finitely many solutions, unless it has an ``obvious'' reason to have infinitely many (reduces to a linear or quadratic equation by a variable change).  

\begin{sloppypar}
Precisely speaking, Siegel proved that an irreducible  equation ${P(x,y)=0}$  (where ${P(x,y)\in \Q[x,y]}$) has at most finitely many solutions  ${x,y\in \Z}$    if one of the following conditions is satisfied:
\end{sloppypar}

\begin{itemize}
\item
the genus of the plane curve ${P(x,y)=0}$ is at least~$1$, or

\item
this curve has at least~$3$ points at infinity. 

\end{itemize}

More generally, let~$\bar C$ be an absolutely irreducible projective curve defined over a number field~$K$ and let~$C$ be an affine subset of~$\bar C$ embedded into the affine space~$\A^\nu$. Further, let~$S$ be a finite set of absolute values of~$K$, including all archimedean absolute values, and let~$\OO_S$ be the ring of $S$-integers of~$K$. Again, Siegel's theorem (in the more general form due to Mahler and Lang) asserts that~$C$ has at most finitely many points in $\A^\nu(\OO_S)$  if ${\genus(\bar C)\ge 1}$ or if ${\left|\bar C\setminus C\right|\ge 3}$. 

Of course, one should mention the celebrated result of Faltings, who proved that the set of rational points on a projective curve of genus~$2$ or higher  is finite. We do not discuss Faltings' work here.

The conventional proof of Siegel's theorem, as in \cite[Chapter~8]{La83} or \cite[Section~D.9]{HS00}, relies on the Theorem of Roth\footnote{At the time of Siegel Roth's theorem was not available, and Siegel had to use a weaker statement.} and heavily depends on the existence of the Jacobian embedding ${\bar C\hookrightarrow J(\bar C)}$, because it exploits high degree étale coverings  of~$\bar C$. 

Recently Corvaja and Zannier~\cite{CZ02} suggested a beautiful new proof, based on the Subspace Theorem rather than the Theorem of Roth, and using projective rather than Jacobian embeddings.

Corvaja and Zannier prove the following theorem.

\begin{theorem}
\label{ttp}
In the above set-up assume that ${\left|\bar C\setminus C\right|\ge 3}$. Then~$C$  has at most finitely many points in $\A^\nu(\OO_S)$.  
\end{theorem}

Siegel's theorem easily follows from Theorem~\ref{ttp}. Indeed, if ${\genus(\bar C)\ge 1}$ then there is an étale covering ${\bar C'\to \bar C}$ of degree~$3$. It induces the covering of affine curves ${C'\to C}$, and we have ${\left|\bar C'\setminus C'\right|\ge 3}$. 

By the Chevalley-Weil principle, the set $\bar C(K)$ is covered by $\bar C'(K')$, where~$K'$ is a number field. Theorem~\ref{ttp} implies that the set of $\OO_{S'}$-integral points on~$C'$ is finite (where~$S'$ is the extension of~$S$ to~$K'$). Hence so is the set of $S$-integral points on~$C$. 

Existence of the covering ${\bar C'\to \bar C}$ of degree~$3$  is the only point in the new proof of Siegel's theorem which appeals to 
the Jacobian embedding: as we shall see, the proof of Theorem~\ref{ttp} is free of Jacobians.

\paragraph*{\textsc{Proof of Theorem~\ref{ttp}}}
Write ${\bar C\setminus C=\{Q_1, \ldots, Q_r\}}$, where, by the assumption, ${r\ge 3}$. Extending the field~$K$, we may assume that each of the points ${Q_1, \ldots, Q_r}$  is defined over~$K$. Further, let ${D=Q_1+ \cdots + Q_r}$ be the ``divisor at infinity''.

Let~$n$ be a (big) positive integer, to be specified later. By the Riemann-Roch theorem, the dimension ${\ell=\ell(nD)}$ of the vector space 
$$
\LL=\LL(nD) = \left\{y \in K(C) : (y)+nD\ge 0\right\}
$$ 
is given by
${\ell= nr-O(1)}$.
In particular, for big~$n$ we have ${\ell\sim nr}$. 

Pick a basis ${y_1, \ldots, y_\ell}$ of~$\LL$. Every~$y_j$ is integral over the   ring ${K[\bx]=K[x_1, \ldots, x_\nu]}$, where ${x_1, \ldots, x_\nu}$ are the coordinate functions on the affine curve ${C\subset \A^\nu}$. 
Multiplying each by a suitable non-zero constant, we may assume that they are integral over the ring $\OO_S[\bx]$. It follows that for every $S$-integral point~$P$ we have ${y_j(P) \in \OO_S}$. 

Now assume that there exist infinitely many distinct $S$-integral points ${P_1, P_2, P_3, \ldots}$. Since~$\bar C$ is a projective curve, the set ${\bar C(K_v)}$ is compact in the $v$-adic topology for every~$v$. Hence, replacing the sequence ${(P_i)}$ by a suitable subsequence, we may assume that it converges in $v$-adic topology for every ${v\in S}$, and we denote by~$Q_v$ the corresponding limits. Now we partition our set~$S$ as ${S=S'\cup S''}$, letting~$S'$ consist of ${v\in S}$ such that ${ Q_v\in \bar C\setminus C}$ and~$S''$ of those~$v$ for which ${ Q_v\in C}$.

We wish to estimate ${\left|y_j(P_i)\right|_v}$ for ${i=1,2,\ldots}$  and ${v\in S}$. For ${v\in S''}$ it is obvious that ${\left|y_j(P_i)\right|_v}$ are bounded independently of~$k$. For ${v\in S'}$ fix a local parameter~$t_v$ at~$ Q_v$. Then ${\left|y_j(P_i)\right|_v\ll \left|t_v(P_i)\right|_v^{-n}}$, where here and below implicit constants are independent of~$i$. Thus, for ${\by=(y_1, \ldots, y_\ell)}$ we obtain
$$
\left\|\by(P_i)\right\|_v \ll
\begin{cases} 
\left|t_v(P_i)\right|_v^{-n}, &\text{if $v\in S'$}\\
1, &\text{if $v\in S''$}.
\end{cases}
$$
Since the numbers ${y_j(P_i)}$ are $S$-integers, we obtain
\begin{equation}
\label{eheig}
H(\by(P_i))=\prod_{v\in S}\left\|\by(P_i)\right\|_v\ll \prod_{v\in S'}\left|t_v(P_i)\right|_v^{-n}.
\end{equation}

All this was just a preparation, and now we are coming to the heart of the Corvaja-Zannier argument. Fix ${v\in S'}$. If ${z\in \LL}$ vanishes at~$ Q_v$, then ${|z(P_i)|_v}$ becomes ``very small'' as~$P_i$ approaches~$ Q_v$, which gives rise to $v$-adically small linear form. Since the vector space~$\LL$ contains ``many'' such~$z$, we have many independent $v$-adically small linear forms. This would  allow us to use the Subspace Theorem.

More specifically, elementary linear algebra shows that our space~$\LL$ has a basis\footnote{It would be more correct to write ${z_{1,v}, \ldots, z_{\ell,v}}$, but this would make the notation too heavy.} ${z_1, \ldots, z_\ell}$ satisfying 
$$
\ord_{ Q_v}z_k\ge k-n-1 \qquad (k=1, \ldots, \ell).
$$
Of course, not all of the functions~$z_k$ vanish at~$ Q_v$ (some of them even have a pole at~$ Q_v$) but, ``in average'', they do. Indeed
\begin{equation}
\label{eposit}
\sum_{k=1}^\ell\ord_{ Q_v}z_k \ge \sum_{k=1}^\ell (k-n-1) = \frac12\ell(\ell-2n-1)=:A.
\end{equation}
Since ${\ell\sim rn}$ for large~$n$, and ${r\ge 3}$ by the assumption, we may specify~$n$ to have ${A>0}$.

Express every~$z_k$ as a linear form in~$\by$: 
$$
z_k=L_{k,v}(\by).
$$
This defines independent linear forms ${L_{1,v}, \ldots, L_{\ell,v}}$ for ${v\in S'}$. For ${v\in S''}$ we simply put ${L_{k,v}(\by)=y_k}$.

We wish to estimate ${\left|L_{k,v}(\by(P_i))\right|_v}$ for all~$k$ and~$v$.  For ${v\in S''}$ we again have 
$$
\left|L_{k,v}(\by(P_i))\right|_v=\left|y_k(P_i)\right|_v\ll 1,
$$ 
and for ${v\in S'}$ we have 
$$
\left|L_{k,v}(\by(P_i))\right|_v = \left|z_k(P_i)\right|_v\ll \left|t_v(P_i)\right|_v^{\ord_{Q_v}z_k}.
$$
Putting this together, we obtain
$$
\prod_{v\in S}\prod_{k=1}^\ell\left|L_{k,v}(\by(P_i))\right|_v\ll \prod_{v\in S'}\left|t_v(P_i)\right|_v^{\sum_{k=1}^\ell \ord_{Q_v}z_k}\le \prod_{v\in S'}\left|t_v(P_i)\right|_v^A,
$$
where ${A>0}$ is defined in~(\ref{eposit}). Combining this with~(\ref{eheig}), we obtain 
$$
\prod_{v\in S}\prod_{k=1}^\ell\left|L_{k,v}(\by(P_i))\right|_v\ll H(\by(P_i))^{-\eps}
$$
with ${\eps=A/n}$. 

Now apply the Subspace Theorem in the form of Theorem~\ref{tssnf}. We obtain that there exist finitely many non-zero functions ${u_1, \ldots, u_s}$ from~$\LL$ such that every~$P_i$ is a zero of one of~$u_j$. It follows that among the points~$P_i$  only finitely many are distinct, which contradicts the original assumption about the existence of an infinite sequence of distinct $S$-integral points. The theorem is proved. \qed

\medskip

Since this argument does not use Jacobians, one may expect to extend to higher dimensions. This is discussed in Subsection~\ref{sssur}. Another useful aspect of the new proof of Siegel's theorem is that it allows, in many cases, to obtain good quantitative bounds for the number of integral points. This direction is exploited, in particular, in~\cite{CZ03}.

\subsection{Integral Points on Surfaces}
\label{sssur}
It is widely believed that an affine (respectively, projective) variety~$V$ of general type cannot have many integral (respectively, rational) points. Of course, one cannot have here ultimate finiteness, but it is expected that integral (or rational) points are not Zariski dense\footnote{Recall that a subset of an algebraic variety is \emph{not Zariski dense} if it lies on a proper closed subvariety.} on~$V$.  Faltings~\cite{Fa91} did the case when~$V$ is a subvariety of an abelian variety, and Vojta extended his result to subvarieties of semiabelian varieties, but very little is known for general~$V$.

Since the argument of Corvaja and Zannier does not use Jacobians, it is very likely to  extend to certain surfaces and varieties of higher dimension, the assumption \textit{there exists at least~$3$ points at infinity} being replaced by something like \textit{the divisor at infinity is ``sufficiently reducible''}. Vojta~\cite{Vo87,Vo96}  used the Subspace Theorem to show that integral points on an irreducible affine variety of dimension~$d$ are not Zariski dense if the divisor at infinity has at least ${d+\rho+1}$ components, where~$\rho$ is the rank of the Néron-Severi group (see also~\cite{NW02}). 

In the article~\cite{CZ04} Corvaja and Zannier applied their argument to integral points on surfaces. Let~$\bar X$ be a non-singular projective surface and ${X\subset \A^\nu}$ a non-empty affine subset of~$\bar X$. We let ${C_1, \ldots, C_r}$ be the irreducible components of ${\bar X\setminus X}$ and we may define the ``divisor at infinity'' ${D=C_1+\cdots+C_r}$. Corvaja and Zannier, however, use the divisor
$$
D=a_1C_1+\cdots+a_rC_r
$$
with some positive integers ${a_1, \ldots, a_r}$ (``weights''). This approach is much more flexible, because the weights can be chosen in a certain ``optimal'' way. 

Recall that in the case of curves we could apply the Subspace Theorem because for every point at infinity~$Q$ and for a sufficiently large~$n$ we found a basis ${z_1, \ldots, z_\ell}$ of the space $\LL(nD)$ such that 
$$
\sum_{j=1}^\ell\ord_Q(z_j)>0.
$$
Similarly, in the surface case, we must find, for every curve~$C_i$ and for a sufficiently large~$n$,   a basis ${z_1, \ldots, z_\ell}$ of the space ${H^0(\bar X, nD)}$ such that 
$$
\sum_{j=1}^\ell\ord_{C_i}(z_j)>0.
$$
We want to express this property in terms of the divisor~$D$. In the subsequent paragraph we  write~$C$ for~$C_i$ and~$a$ for~$a_i$.

Consider the filtration of the space $H^0(\bar X,nD)$
\begin{equation}
\label{efilt}
H^0(\bar X,nD)\supseteq H^0(\bar X,nD-C)\supseteq H^0(\bar X,nD-2C)\supseteq\ldots ,
\end{equation}
and let ${z_1, \ldots, z_\ell}$ be a basis of this filtration\footnote{A basis of a filtration ${W_0\supseteq W_1\supseteq W_2\supseteq\ldots}$ of  vector spaces is, by definition, a basis of~$W_0$ which contains a basis of every~$W_i$.}. For this basis we have
\begin{align*}
\sum_{j=1}^\ell\ord_{C}(z_j) &= \sum_{k=0}^\infty (k-an)\Bigl(h^0\bigl(nD-kC\bigr) -h^0\bigl(nD-(k+1)C\bigr)\Bigr)\\ &= -anh^0(nD)+ \sum_{k=0}^\infty h^0(nD-kC) 
\end{align*}
(of course, the infinite sums have only finitely many non-zero terms).

Thus, the basic condition to be satisfied  is that  the inequalities
\begin{equation}
\label{efract}
\frac{\sum_{k=0}^\infty h^0(nD-kC_i)}{nh^0(nD)} >a_i \qquad (i=1, \ldots, r)
\end{equation}
hold for a certain~$n$. 

\begin{rtheorem}{Corvaja, Zannier}
\label{tczs}
Let~$\bar X$ be a non-singular projective surface defined over a number field~$K$ and let ${X\subset \A^\nu}$ be a non-empty affine subset of~$\bar X$. Let ${C_1, \ldots, C_r}$  be effective divisors\footnote{We do not assume the divisors ${C_1, \ldots, C_r}$ irreducible.} supported at ${\bar X\setminus X}$. Assume that ${C_1, \ldots, C_r}$ intersect properly (that is, no~$2$ of them have a common component and no~$3$ of them have a common point). Further, assume that for some choice of positive integers ${a_1, \ldots, a_r}$ the~$r$ inequalities~(\ref{efract}) (with ${D=a_1C_1+\cdots+a_rC_r}$) hold for certain~$n$. 
Then for any finite set ${S\subset M_K}$ the set ${X\cap \A^\nu(\OO_S)}$ of $S$-integral points on~$X$ is not Zariski dense. 
\end{rtheorem}

\paragraph*{\textsc{Proof}} 
It is quite analogous to the proof of Theorem~\ref{ttp}. We may assume that every~$C_i$ is defined over~$K$. Let~$n$ be such that the inequalities~(\ref{efract}) hold. As we have seen above, this implies existence of a positive~$B$ such that 
$$
\sum_{k=1}^\ell\ord_{C}z_k \ge B, 
$$
where~$C$ is any of ${C_1, \ldots, C_r}$ and ${z_1, \ldots, z_\ell}$ is a basis of the filtration~(\ref{efilt}).

To prove the theorem, it suffices to show that every infinite sequence of $S$-integral points has a subsequence contained on a curve defined over~$K$. Indeed, since there is only countably many $K$-curves, a Zariski-dense set contains a sequence with finitely many elements on every $K$-curve. 

Thus, let ${P_1,P_2,P_3 \ldots}$ be sequence of $S$-integral points. Replacing it by a subsequence, we may assume that it $v$-adically converges for every ${v\in S}$, and denote the limit by~$Q_v$. Now we have~$3$ cases:
either ${Q_v \in X}$ or~$Q_v$ belongs exactly one of the~$C_i$ (call it~$C_v$), or it belongs to exactly two of them (call them $C_v$ and $C_v'$). (By the assumption,~$Q_v$ cannot belong to three or more of~$C_i$.) Let~$S_0$,~$S_1$ and~$S_2$ be the corresponding subsets of~$S$. 

Fix a basis ${y_1, \ldots, y_\ell}$ of the space ${H^0(\bar X,nD)}$. We may assume that ${y_j(P)\in \OO_S}$ for any $S$-integral point~$P$.

Now, for each ${v\in S}$ we shall define a new basis ${z_1=z_{1,v}, \ldots, z_\ell=z_{\ell,v}}$ of the same space, and we let ${L_{1,v}, \ldots, L_{\ell, v}}$ be the linear forms such that 
${z_k=L_{k,v}(\by)}$. Then we shall apply the Subspace Theorem to these forms evaluated at  $\by(P_i)$.

If ${v\in S_0}$ then, as in the proof of Theorem~\ref{ttp}, we define the $z$-basis just putting ${z_j=y_j}$. We have plainly 
\begin{align}
\label{esny}
\|\by(P_i)\|_v &\ll 1, \\
\label{esnz}
\prod_{k=1}^\ell\left|L_{k,v}(\by(P_i))\right|_v&\ll 1.
\end{align}

Next, assume that ${v \in S_1}$ and let ${z_1, \ldots, z_\ell}$ be a basis of the filtration~(\ref{efilt}) with ${C=C_v}$. If~$t_v$ is a local parameter of~$C_v$ near~$Q_v$ then for any function~$u$ regular on~$X$ the function ${t_v^{-\ord_{C_v}u}u}$ is regular in a neighborhood of~$Q_v$. It follows 
$$
|u(P_i)|_v \ll |t_v(P_i)|_v^{\ord_{C_v}u} \qquad (i=1,2, \ldots).
$$
Applying this with ${u=y_1, \ldots, y_\ell}$ and with ${u=z_1, \ldots, z_\ell}$, we find that 
\begin{alignat}2
\label{esoy}
\|\by(P_i)\|_v &\ll |t_v(P_i)|_v^{\min_{1\le j\le \ell} \ord_{C_v}y_j}&&\le |t_v(P_i)|_v^{-An},\\
\label{esoz}
\prod_{k=1}^\ell\left|L_{k,v}(\by(P_i))\right|_v&\ll |t_v(P_i)|_v^{\sum_{k=1}^\ell\ord_{C_v}z_k} &&\le |t_v(P_i)|_v^B,
\end{alignat}
where ${A=\max\{a_1, \ldots, a_r\}}$ and ${B>0}$ is defined in the beginning of the proof.

Finally, assume that ${v\in S_2}$. In this case  Corvaja and Zannier use the following nice elementary lemma.

\begin{lemma} 
\label{lfilt}
Let
\begin{equation}
\label{efilem}
W=W_0\supseteq W_1\supseteq W_2\supseteq\ldots, \qquad 
W=W_0'\supseteq W_1'\supseteq W_2'\supseteq\ldots
\end{equation} 
be two filtrations of a finitely dimensional vector space~$W$. Then there exists a common basis for the two filtrations (That is, there exists a basis of~$W$ containing bases for every~$W_i$ and for every~$W_i'$.) 
\end{lemma}

{\footnotesize (The proof is by induction in $\dim W$.  Without loss of generality we may assume that~$W_1$ is a hyperplane in~$W$. Put ${W_i''=W_1\cap W_i'}$. By induction, there exists a common basis ${w_1, \ldots, w_{d-1}}$ for the filtrations ${ W_1\supseteq W_2\supseteq\ldots}$ and ${W_1=W_0''\supseteq W_1''\supseteq W_2''\supseteq\ldots}$. Now let~$k$ be the smallest index for which ${W_k'\not\subseteq W_1}$ (the set of such indices  is non-empty because it includes~$0$). Then~$W_i''$ is a hyperplane in~$W_i'$ for ${i\le k}$ and ${W_i'=W_i''}$ for ${i>k}$. Now, picking a ${w_d\in W_k'\setminus W_k''}$, we obtain a basis ${w_1, \ldots, w_{d-1}, w_d}$ of both filtrations~(\ref{efilem}), which proves the lemma.)}

\medskip
Using the lemma, we find a common basis ${z_1, \ldots, z_\ell}$ for both the filtrations~(\ref{efilt}) with ${C=C_v}$ and ${C=C_v'}$. 

Now let~$t_v$ and~$t_v'$ be local parameters near~$Q_v$ at~$C_v$ and~$C_v'$, respectively. Then for any function~$u$ regular on~$X$ the function ${t_v^{-\ord_{C_v}u}(t_v')^{-\ord_{C_v'}u}u}$ is regular in a neighborhood of~$Q_v$, whence 
$$
|u(P_i)|_v \ll |t_v(P_i)|_v^{\ord_{C_v}u}|t_v'(P_i)|_v^{\ord_{C_v'}u} \qquad (i=1,2, \ldots).
$$
Applying this with ${u=y_1, \ldots, y_\ell}$ and with ${u=z_1, \ldots, z_\ell}$, we obtain
\begin{alignat}2
\label{esty}
\|\by(P_i)\|_v &\ll |t_v(P_i)|_v^{\min_{1\le j\le \ell} \ord_{C_v}y_j}|t_v'(P_i)|_v^{\min_{1\le j\le \ell} \ord_{C_v'}y_j}&&\le \left|t_v(P_i)t_v'(P_i)\right|_v^{-An},\\
\label{estz}
\prod_{k=1}^\ell\left|L_{k,v}(\by(P_i))\right|_v&\ll |t_v(P_i)|_v^{\sum_{k=1}^\ell\ord_{C_v\vphantom{C_v'}}z_k} |t_v'(P_i)|_v^{\sum_{k=1}^\ell\ord_{C_v'}z_k}&&\le |t_v(P_i)t_v'(P_i)|_v^B.
\end{alignat}

Combining the inequalities (\ref{esny},\ref{esoy},\ref{esty}) with (\ref{esnz},\ref{esoz},\ref{estz}), we find 
$$
\prod_{v\in S}\prod_{k=1}^\ell\left|L_{k,v}(\by(P_i))\right|_v\ll H(\by(P_i))^{-\eps} \qquad (i=1,2,\ldots)
$$
with ${\eps=B/An}$. Now we complete the proof using the Subspace Theorem in the same manner as we did in the proof of Theorem~\ref{ttp}. \qed

\begin{remark}
Using Vojta's refinement~\cite{Vo89} of the Subspace Theorem, Levin~\cite{Le07} shows that, under the hypothesis of Theorem~\ref{tczs} there exists a (possibly, reducible) affine curve on~$X$,  depending only on~$X$, but independent on~$K$ and~$S$, such that all but finitely many $S$-integral points from~$X$ belong to this curve. (The exceptional finite set may, however, depend on~$K$ and~$S$.) The same is true for the consequences of Theorem~\ref{tczs}: Corollaries~\ref{ccz} and~\ref{cautis} and Theorem~\ref{tlev}. 
\end{remark}

Imposing on our divisors~$C_i$  additional assumption (like ampleness), we can estimate from below  the quantity on the left of~(\ref{efract}) asymptotically   (as ${n\to \infty}$), using the Riemann-Roch theorem on surfaces.  This would express our condition in terms of the intersection numbers of the divisors ${C_1, \ldots, C_r}$ and the weights ${a_1, \ldots, a_r}$. The Riemann-Roch theorem applies through the following lemma, proved in the addendum to this section. 

\begin{sloppypar}
	\begin{lemma}
	\label{lrr}
	Let~$C$ be an ample divisor and~$D$ an effective divisor on a non-singular projective surface, and let~$n$ and~$k$ be positive integers such that ${k\le \alpha n}$, where 
	${\alpha =(D\cdot C)/C^2}$.
	Then 
	\begin{equation}
	\label{equad}
	h^0(nD-kC)\ge \frac12(nD-kC)^2 -O(n).
	\end{equation}
	\end{lemma}
\end{sloppypar}

Let us look closer at this lemma. We have
$$
(nD-kC)^2=D^2n^2-2(D\cdot C) nk+ C^2k^2. 
$$
The quadratic form
$$
q(\xi, \tau)= D^2\xi^2-2(D\cdot C) \xi\tau+ C^2\tau^2
$$
is not positive definite by the Hodge index theorem. Hence the polynomial 
$q(1, \tau)$ has two real roots,~$\gamma$ and~$\gamma'$. They are, obviously, positive, and we assume that ${\gamma\le \gamma'}$. In fact, ${\gamma\le \alpha\le \gamma'}$ because ${\alpha=(\gamma+\gamma')/2}$.

Thus, we have ${q(\xi, \tau)<0}$ if ${\gamma \xi<\tau<\gamma'\xi}$, and ${q(\xi, \tau)\ge0}$ otherwise. In particular,~(\ref{equad}) remains true for ${k\le\gamma'n}$, but it is is uninteresting for ${\gamma n<k<\gamma' n}$ and becomes interesting only for ${k\le \gamma n}$.

Applying the lemma in our situation, we bound the numerator on the left of~(\ref{efract}) as (we write~$C$ instead of~$C_i$)
\begin{align}
\label{etheta}
\sum_{k=0}^\infty h^0(nD-kC) &\ge \sum_{k\le \theta n}\frac12q(n,k) -O(n^2) \\
&= \left(\frac12 \theta D^2 -\frac12\theta^2(D\cdot C) +\frac16\theta^3C^2\right)n^3 -O(n^2), \nonumber
\end{align} 
where~$\theta$ is any real number satisfying ${0\le \theta \le \gamma'}$. Also, the Riemann-Roch theorem gives for the denominator in~(\ref{efract}) the asymptotics
$$
nh^0(nD)= \frac12n^3D^2 + O(n^2).
$$
Hence the left-hand side of~(\ref{efract}) is bounded from below by ${F(\theta)+ O(1/n)}$, where 
$$
F(\theta)=\theta\left(1-\theta \frac{D\cdot C}{D^2} + \frac13 \theta^2\frac{C^2}{D^2}\right) .
$$
It remains to select the parameter~$\theta$ in the optimal way.

The estimate in~(\ref{etheta})  is best possible if the sum on the right of~(\ref{etheta}) contains all positive terms $q(n,k)$ and no negative terms. It follows that the optimal choice is ${\theta =\gamma}$. We obtain the following consequence. 

\begin{rcorollary}{Corvaja, Zannier}
\label{ccz}
Let~$\bar X$ be a non-singular projective surface defined over a number field~$K$ and let ${X\subset \A^\nu}$ be a non-empty affine subset of~$\bar X$. Let ${C_1, \ldots, C_r}$  be properly intersecting effective ample divisors supported at ${\bar X\setminus X}$.  Further, assume that for some choice of positive integers ${a_1, \ldots, a_r}$ the~$r$ inequalities
\begin{equation}
\label{ecz}
\gamma_i\left(1-\gamma_i \frac{D\cdot C_i}{D^2} + \frac13 \gamma_i^2\frac{C_i^2}{D^2}\right) >a_i  \qquad (i=1, \ldots r)
\end{equation}
hold, where ${D=a_1C_1+\cdots +a_rC_r}$ and where~$\gamma_i$ is the smallest positive root of the polynomial ${D^2-2(D\cdot C_i)T+ C_i^2T^2}$. 
Then for any finite set ${S\subset M_K}$ the  $S$-integral points are not Zariski dense on~$X$. 
\end{rcorollary}

By choosing suitable weights, Corvaja and Zannier showed that integral points are not Zariski dense if satisfy some condition; for instance, if ${r\ge 4}$ and  the intersection matrix of ${C_1, \ldots, C_r}$ is of rank~$1$. 

Autissier~\cite{Au07} suggested to take ${\theta = \beta/2}$, where 
${\beta = D^2/(D\cdot C)}$.
(Notice that ${\beta/2 <\gamma}$ and ${\gamma\approx \beta/2}$ when~$\gamma'$ is very large.) Since 
$$
F\left(\frac\beta2\right)= \frac\beta2\left(1-\frac\beta2 \frac{D\cdot C}{D^2} + \frac{\beta^2}{12}\frac{C^2}{D^2}\right)=
\frac 14 \frac{D^2}{D\cdot C}\left(1+ \frac16 \frac{D^2C^2}{(D\cdot C)^2}\right),
$$
we obtain the following result.

\begin{rcorollary}{Autissier}
\label{cautis}
In the set-up of Corollary~\ref{ccz}, assume that for some choice of positive integers ${a_1, \ldots, a_r}$ the~$r$ inequalities
\begin{equation}
\label{eaut}
\frac{D^2}{D\cdot C_i}\left(1+ \frac16 \frac{D^2C_i^2}{(D\cdot C_i)^2}\right) >4a_i  \qquad (i=1, \ldots, r)
\end{equation}
hold. 
Then for any finite set ${S\subset M_K}$ the  $S$-integral points are not Zariski dense on~$X$. 
\end{rcorollary}

This result is formally weaker, than Corollary~\ref{ccz}, but it is more practical, because inequality~(\ref{eaut}) is much easier to handle, than~(\ref{ecz}).

Levin~\cite{Le07}, and, independently, Autissier~\cite{Au07} observed that a ``nearly optimal'' choice  of the weights ${a_1, \ldots, a_r}$ implies that~$4$ ample divisors at infinity would suffice. More precisely, they prove the following. 

\begin{rtheorem}{Levin, Autissier}
\label{tlev}
Let~$\bar X$ be a non-singular projective surface defined over a number field~$K$ and let ${X\in \A^\nu}$ be a non-empty affine subset of~$\bar X$. Let ${C_1, \ldots, C_r}$  be properly intersecting effective ample divisors supported at ${\bar X\setminus X}$. Assume that ${r\ge 4}$. Then for any finite set ${S\subset M_K}$ the $S$-integral points on~$X$ are not Zariski dense. 
\end{rtheorem}

\begin{remark}
In Theorem~\ref{tlev} one can relax the assumption that the divisors~$C_i$ are ample (see \cite[Theorem~11.5A]{Le07}), but one cannot just assume that~$C_i$ are effective and intersect properly. As an example take ${\bar X=\P^1\times \P^1}$ and ${X=\G_m\times\G_m}$, where $\G_m$ is obtained by removing the $0$-point and the $\infty$-point from~$\P^1$. Then ${\bar X\setminus X}$ consists of~$4$ curves. The map ${(x,y)\to (x,x^{-1},y,y^{-1})}$ defines an affine embedding ${X\to \A^4}$, and the set of $S$-integral points with respect to this embedding is $\OO_S^\times\times \OO_S^\times$, which is Zariski-dense in general. 
\end{remark}

To prove Theorem~\ref{tlev} we need one more elementary lemma.

\begin{lemma}
\label{laut}
Let ${M=\left[\mu_{ij}\right]_{1\le i,j\le r}}$ be a symmetric ${r\times r}$-matrix with positive real entries. Consider the linear forms
$$
L_i(\bx)=\mu_{i1}x_1+\cdots+\mu_{ir}x_r \qquad (i=1, \ldots, r)
$$
and the quadratic form ${Q(\bx)=\bx^{t}M\bx}$. Then for any ${\eps>0}$ there exist positive integers ${a_1, \ldots, a_r}$ such that 
\begin{equation}
\label{eautis}
(1-\eps) Q(\ba)< ra_iL_i(\ba) < (1+\eps)Q(\ba)   \qquad (i=1, \ldots, r),
\end{equation}
where ${\ba =(a_1, \ldots, a_r)}$. 
\end{lemma}

\paragraph*{\textsc{Proof}}
We follow the elegant argument of Autissier \cite[Proposition~2.3]{Au07}. Notice that 
$$
Q(\bx) = x_1L_1(\bx) + \cdots + x_rL_r(\bx). 
$$
Hence we have to find a point~$\ba$ with  positive \emph{integral} coordinates such that the~$r$ numbers $a_iL_i(\ba)$ are \emph{approximately} equal. We first find a point with positive \emph{real} coordinates where these numbers are \emph{exactly} equal. 

Let~$\Delta$ be the simplex 
\begin{equation}
\label{esimple}
x_1+\cdots+x_r=1, \qquad 0\le x_i\le 1 \qquad (i=1,\ldots, r).
\end{equation}
Consider the map 
${\Delta\to \Delta}$ defined by 
$$
\bx \mapsto \left(L_1(\bx)^{-1}, \ldots, L_r(\bx)^{-1}\right)\left(\sum_{i=1}^rL_i(\bx)^{-1}\right)^{-1}.
$$
The map is well-defined because the entries of our matrix~$M$ are positive numbers. By the Brower theorem, our map has a fixed point ${\ba \in \Delta}$. For this point we have 
$$
a_1L_1(\ba)=\ldots =a_rL_r(\ba).
$$
Since none of the $L_i(\ba)$ vanishes, none of the ~$a_i$ does; in other words, the real numbers ${a_1, \ldots, a_r}$ are strictly positive. Replacing each by a suitable rational approximation, we obtain positive rational numbers ${a_1, \ldots, a_r}$ satisfying~(\ref{eautis}). Multiplying them by the common denominator, we arrive to the desired integers ${a_1, \ldots, a_r}$. \qed

\medskip

\paragraph*{\textsc{Proof of Theorem~\ref{tlev}}}
First of all, remark that the term 
\begin{equation}
\label{eterm}
\frac16 \frac{D^2C_i^2}{(D\cdot C_i)^2},
\end{equation}
occurring in~(\ref{eaut}), is bounded from below, uniformly in~$\ba$, by a positive constant.   Indeed,~(\ref{eterm}) defines a homogeneous positive real function on non-zero vectors   ${\ba\in\left(\Z_{\ge0}\right)^r}$. But, since it is a quotient of quadratic forms with positive coefficients, it extends to a positive real continuous function  on  the non-zero vectors of ${\left(\R_{\ge0}\right)^r}$ By homogeneity, it suffices to consider this function on  the compact~$\Delta$ defined by~(\ref{esimple}), where it is bounded away from~$0$.

Thus, to ensure~(\ref{eaut}), we must find positive integers ${a_1, \ldots, a_r}$ such that for some ${\eps>0}$ the inequalities 
$$
D^2\left(1+ \eps\right) >4a_i(D\cdot C_i)  \qquad (i=1, \ldots, r)
$$
hold. Applying Lemma~\ref{laut} to the intersection matrix of ${C_1, \ldots, C_r}$, we find ${a_1, \ldots, a_r}$ such that
$$
D^2\left(1+ \eps\right) >ra_i(D\cdot C_i)  \qquad (i=1, \ldots, r).
$$
Since ${r\ge 4}$, we are done. \qed

\medskip

In his fundamental article~\cite{Le07} Levin extends Theorem~\ref{tlev} to varieties of arbitrary dimension, without assuming proper intersection. One difficulty he has to overcome is that Lemma~\ref{lfilt} is no longer true for three or more filtrations. 

Levin gives a thorough analysis of the argument of Corvaja and Zannier and, probably,  reaches its ``natural limitations''. In addition, he accompanies every Diophantine result with  an analogous statement about holomorphic maps, in accordance with Vojta's philosophy.

For more Diophantine applications of the Subspace Theorem see~\cite{CZ04a,CZ06}.

\subsection*{Addendum: Proof of Lemma~\ref{lrr}}
We deduce Lemma~\ref{lrr} from the Theorem of Riemann-Roch and the following proposition.

\begin{proposition}
\label{pbcd}
Let~$B$,~$C$ and~$D$  be  divisors on a non-singular projective surface~$X$. Assume that~$C$ is very ample, that~$D$ is effective and that ${C^2\le C\cdot D}$.  Then
$$
h^0(B-D+C) \le B\cdot C+h^0(B)+1.
$$
\end{proposition}

\paragraph*{\textsc{Proof}}
By the Theorem of Bertini we may assume that~$C$ is an irreducible smooth curve.  The exact sequence of sheaves 
$$
0\to \OO_X(B-D)\to \OO_X(B-D+C) \to \OO_X(B-D+C)\vert_C\to 0,
$$
implies the exact sequence of cohomologies
$$
0\to H^0(X, B-D)\to H^0(X,B-D+C)\to H^0(C, \Delta)\to \ldots,
$$
where~$\Delta$ is the divisor ${(B-D+C)\vert_C}$ on~$C$. It follows that
\begin{equation}
\label{eboun}
h^0(X,B-D+C)\le h^0(X,B-D) + h^0\left(C,\Delta\right).
\end{equation}
We have
${\deg \Delta= (B-D+C)\cdot C \le B\cdot C}$
because ${C^2\le C\cdot D}$. It remains to observe that 
${h^0\left(C,\Delta\right)\le \deg \Delta+1}$ and that  
${h^0(X,B-D)\le h^0(X,B)}$, 
because~$D$ is effective.  \qed

\medskip

\paragraph*{\textsc{Proof of Lemma~\ref{lrr}}}
By the theorem of Riemann-Roch,
$$
h^0(nD-kC) \ge \frac12(nD-kC)^2 -\frac12\bigl((nD-kC)\cdot K\bigr) - h^0(K-nD+kC)+ O(1),
$$
where~$K$ is the canonical divisor. 
Since ${(nD-kC)\cdot K =O(n)}$, we have to prove that ${h^0(K-nD+kC)=O(n)}$. 

We may assume~$k$ so large that~$kC$ is very ample. Applying Proposition~\ref{pbcd} with ${B=K}$ and with~$kC$,~$nD$ instead of~$C$ and~$D$,
the condition ${(kC)^2\le kC\cdot nD}$ being assured by the assumption ${k\le \alpha n}$, we find 
$$
h^0(K-nD+kC)\le k(K\cdot C)+h^0(K)+1=O(n),
$$
as wanted.\qed

\medskip

I thank Ivan Cheltsov for explanations concerning this lemma.

\section{Conclusion}
As it was indicated in the introduction, the recent remarkable applications of the Subspace Theorem are not limited to the results discussed above. Without any claim for exhaustiveness, let me just quote several more works that I personally find attractive. 

Mahler~\cite{Ma57} showed, using the theorem of Ridout, that if~$\alpha$ is a positive  rational number, but not integer, and ${0<\theta<1}$ then the inequality ${\left|\alpha^n-m\right|\le \theta^n}$ has finitely many solutions in positive integers~$n$ and~$m$. He asked for which irrational algebraic numbers a similar statement is true, observing that is is false, for instance, if ${\alpha=\left(1+\sqrt5\right)/2}$, and, more generally, if~$\alpha$ is a Pisot number\footnote{See footnote~\ref{fpis} on page~\pageref{fpis}.}. Corvaja and Zannier
answered this question, showing that the corresponding statement is true for all irrational algebraic numbers except the roots of Pisot numbers (and for the latter it is obviously false). 

In the same article they answered a question of Mendès France~\cite{Me93} on the period length of the periodic continued fraction for~$\alpha^n$, where~$\alpha$ is a quadratic irrationality. Corvaja and Zannier showed that the period tends to infinity with~$n$ unless~$\alpha$ is a square root of a rational number or a unit. See also~\cite{BL05,CZ05,Sc06}. 

Corvaja and Zannier~\cite{CZ02b} gave a complete answer to Pisot's question on when  the quotient ${u(n)/v(n)}$ of two power sums (and, more generally, of two linear recurrences) is infinitely often an integer. By the way, this is one of the rare cases when the authors managed to overcome the difficulty stemming from the absence of the ``dominant root'' (see the end of Section~\ref{spowsum}). 

Corvaja and Zannier~\cite{CZ03a} and, independently, Hernández and Luca~\cite{HL03} proved that  ${(ab+1)(ac+1)(bc+1)}$ cannot have only small prime divisors, confirming a conjecture of Gy\H ory, Sárközy and Stewart~\cite{GSS96}. See~\cite{BL04} for a quantitative version of this result. 

Bugeaud, Corvaja and Zannier proved~\cite{BCZ03} that ${a^n-1}$ and ${b^n-1}$ cannot have a large common divisor. This was extended by Corvaja and Zannier~\cite{CZ03a,CZ05a}.

Corvaja, Rudnick and Zannier~\cite{CRZ04} showed that (with obvious exceptions) the multiplicative order of an  integral matrix $\mod N$ grows quicker than $\log N$ as ${N\to \infty}$. This result is essentially best possible.

And there are numerous  contributions that I failed to mention, because of lack of space or time or because of my ignorance.

{\footnotesize\paragraph*{Acknowledgments} Mayeul Bacquelin explained me the work of Adamczewski and Bugeaud. Umberto Zannier was very helpful and patient when clarifying me various aspects of his work with Corvaja. I also had useful correspondence and/or discussions with Pascal Autissier, Boris Adamczewski, Yann Bugeaud, Ivan Cheltsov, Pietro Corvaja, Aaron Levin and  Hans Peter Schlickewei. Many colleagues, including 
Boris Adamczewski, Yann Bugeaud, Ivan Cheltsov, Pietro Corvaja, Viviane le Dret, Marina Prokhorova and Umberto Zannier,  read the manuscript and detected a number of inaccuracies. I am happy to thank them all.

In preparation of this text, I benefited a lot from Zannier's excellent notes~\cite{Za03}, and I strongly recommend them to anybody wishing to learn more on the Diophantine aspect of the Subspace Theorem.

My deepest gratitude goes to Elina Wojciechowska, for her constant encouragement during my work on this article. }

\end{document}